\documentclass[a4paper]{amsproc}
\usepackage{amssymb}
\usepackage{mathrsfs}
\usepackage{amsfonts}
\usepackage{amsthm}
\usepackage{amsmath,amscd}
\usepackage{amsxtra}     % Use various AMS packages
\usepackage{bm}
\usepackage[all,cmtip]{xy}
\usepackage{epsfig}
\usepackage{verbatim}
\usepackage{color}
\usepackage{enumerate}
\usepackage{hyperref}
\usepackage{tikz}
\usetikzlibrary{arrows,matrix,shapes,trees}
\usepackage{extarrows}
\usepackage{enumitem}

\theoremstyle{plain}
 \newtheorem{thm}{Theorem}[section]
 \newtheorem{prop}[thm]{Proposition}
 \newtheorem{lem}[thm]{Lemma}

\theoremstyle{definition}
 \newtheorem{ex}{Example}[section]
 \newtheorem{defn}{Definition}[section]
\theoremstyle{remark}
 \newtheorem{rmk}{Remark}[section]
 \numberwithin{equation}{section}

\newcommand{\N}{{\mathbb N}}
\newcommand{\Q}{{\mathbb Q}}
\newcommand{\Z}{{\mathbb Z}}

\newcommand{\Gm}{\mathbb{G}_{\mr{m}}}

\newcommand{\Gml}{\mathbb{G}_{\mr{m,log}}}

\newcommand{\mr}{\mathrm}
\newcommand{\mc}{\mathcal}

\newcommand{\Spec}{\mathop{\mr{Spec}}}

\makeatletter
\@namedef{subjclassname@2020}{%
  \textup{2020} Mathematics Subject Classification}
\makeatother

\title[The higher direct images of locally constant group schemes]{The higher direct images of locally constant group schemes from the Kummer log flat topology to the classical flat topology}



\author[Heer Zhao]{\bfseries Heer Zhao}

\address{
    Heer Zhao, 
    Fakult\"at f\"ur Mathematik, 
    Universit\"at Duisburg-Essen, 
    Essen 45117, 
    Germany, 
    	heer.zhao@uni-due.de}

\begin{document}

%{\begin{flushleft}\baselineskip9pt\scriptsize
%PUBLICATIONS DE L'INSTITUT MATH\'EMATIQUE\newline
%Nouvelle s\'erie, tome 91(105) (2012), od--do \hfill DOI:
%\end{flushleft}}
\vspace{18mm} \setcounter{page}{1} \thispagestyle{empty}

\begin{abstract}
Let $S$ be an fs log scheme, and let $F$ be a group scheme over the underlying scheme which is \'etale locally representable by (1) a finite dimensional $\Q$-vector space, or (2) a finite rank free abelian group, or (3) a finite abelian group. We give a full description of all the higher direct images of $F$ from the Kummer log flat site to the classical flat site. In particular, we show that: in case (1) the higher direct images of $F$ vanish; and in case (2) the first higher direct image of $F$ vanishes and the $n$-th ($n>1$) higher direct image of $F$ is isomorphic to the $(n-1)$-th higher direct image of $F\otimes_{\Z}\Q/\Z$. In the end, we make some computations when the base is a standard log trait or a Dedekind scheme endowed with the log structure associated to a finite set of closed points.
\end{abstract}

\maketitle

\section{Introduction}\label{sec1}
Let $S$ be an fs log scheme, $(\mr{fs}/S)$ the category of fs log schemes over $S$. We endow $(\mr{fs}/S)$ with the Kummer log flat topology (resp. the classical flat topology, resp. the classical \'etale topology), see \cite[\S 2]{kat2} (resp. \cite[\S 4]{kat2}, resp. \cite[\S 2.4]{ill1}), and denote the resulting site by $(\mr{fs}/S)_{\mr{kfl}}$ (resp. $(\mr{fs}/S)_{\mr{fl}}$, resp. $(\mr{fs}/S)_{\mr{\acute{e}t}})$\footnote{In \cite{kat2} the sites $(\mr{fs}/S)_{\mr{kfl}}$ and $(\mr{fs}/S)_{\mr{fl}}$ are denoted as $S_{\mr{fl}}^{\mr{log}}$ and $S_{\mr{fl}}^{\mr{cl}}$ respectively. Our notation here is analogous to that of the Kummer log \'etale site from \cite[\S 5.3]{k-k-n4}}. We have a canonical map
\[\varepsilon_{\mr{fl}}:(\mr{fs}/S)_{\mr{kfl}}\to (\mr{fs}/S)_{\mr{fl}}\]
of sites. To understand the cohomology of a sheaf of abelian groups $F$ on the site $(\mr{fs}/S)_{\mr{kfl}}$, one needs to understand the higher direct images $R^i\varepsilon_{\mr{fl}*}F$. 

The first higher direct image $R^1\varepsilon_{\mr{fl}*}F$ has been determined by Kato when $S$ is locally noetherian and $F$ is either a finite flat (commutative) group scheme or a smooth affine (commutative) group scheme over the underlying scheme of $S$, see \cite[Thm. 4.1]{kat2} or \cite[Thm. 3.12]{niz1}. Kato's theorem about $R^1\varepsilon_{\mr{fl}*}F$ has been generalized to quasi-projective smooth (commutative) group schemes by the author, see \cite[Thm. 3.14]{zha5}.

Assume that the underlying scheme of $S$ is locally noetherian. The second higher direct image $R^2\varepsilon_{\mr{fl}*}F$ has been described in \cite[Thm. 3.23]{zha5} (resp. \cite[Thm. 1.2]{zha6}, resp. \cite[Thm. 1.3]{zha6}, resp. \cite[Thm. 1.4]{zha6}) when $F$ is representable by a torus (resp. a smooth affine commutative group scheme, resp. a finite flat commutative group scheme, resp. an extension of an abelian scheme by a torus) over the underlying scheme of $S$. When $F$ is representable by a smooth quasi-projective commutative group scheme, the higher direct images are always torsion by \cite[Thm. 1.1]{zha6}.

In this article, we investigate $R^i\varepsilon_{\mr{fl}*}F$ for all $i>0$ in Section \ref{sec2} when $F$ is representable by a group scheme which is \'etale locally isomorphic either to a finite dimension $\Q$-vector space, or a finite rank free abelian group, or a finite abelian group. The main results are the following three theorems.

\begin{thm}\label{1.1}(See also Theorem \ref{2.9})
Let $S$ be a locally noetherian fs log scheme. Let $l$ be a prime number, $U$ the open locus on $S$ where $l$ is invertible, and $j:U\hookrightarrow S$ the corresponding strict open immersion. Let $F$ be a finite \'etale group scheme over the underlying scheme of $S$, and we endow it with the induced log structure from $S$. Assume that $F$ is killed by an $l$-power, then we have
\[R^i\varepsilon_{\mr{fl}*}F\cong j_{\mr{fl!}}((j_{\mr{kfl}}^{-1}F)(-i)\otimes_{\Z}\bigwedge^i(\Gml/\Gm)_{U_{\mr{kfl}}})\]
for $i\geq1$, where $j_{\mr{kfl}}:(\mr{fs}/U)_{\mr{kfl}}\to (\mr{fs}/S)_{\mr{kfl}}$ (resp. $j_{\mr{fl}}:(\mr{fs}/U)_{\mr{fl}}\to (\mr{fs}/S)_{\mr{fl}}$) is the morphism on the Kummer log flat sites (resp. the classical flat sites) induced by $j$.
\end{thm}

\begin{thm}\label{1.2}(See also Theorem \ref{2.10})
Let $S$ be an fs log scheme. Let $F$ be a group scheme over the underlying scheme of $S$ which is \'etale locally representable by a finite dimensional $\Q$-vector space. Then we have $R^i\varepsilon_{\mr{fl}*}F=0$ for $i\geq 1$.
\end{thm}

\begin{thm}\label{1.3}(See also Theorem \ref{2.11})
Let $S$ be a locally noetherian fs log scheme, and $F$ a group scheme over the underlying scheme of $S$ which is \'etale locally isomorphic to a finite rank free abelian group. Then we have the following.
\begin{enumerate}
\item $R^1\varepsilon_{\mr{fl}*}F=0$.
\item Let $i>1$. For each prime number $l$, let $U_l$ be the locus on $S$ on which $l$ is invertible and ${}_lj:U_l\hookrightarrow S$ the corresponding strict open immersion. Then 
\begin{align*}
R^i\varepsilon_{\mr{fl}*}F&\cong\bigoplus_{l\text{ prime}}R^{i-1}\varepsilon_{\mr{fl}*}(F\otimes_\Z\Q_l/\Z_l) \\
&\cong \bigoplus_{l\text{ prime}}{}_lj_{\mr{fl}!}({}_lj^{-1}F\otimes_\Z\Q_l/\Z_l(-i+1)\otimes_{\Z}\bigwedge^{i-1}(\Gml/\Gm)_{(U_{l})_{\mr{fl}}}).
\end{align*}
\end{enumerate} 
\end{thm}

The proof of Theorem \ref{1.3} is reduced to Theorem \ref{1.1} and Theorem \ref{1.2} via the short exact sequence 
\[0\to F\to F\otimes_{\Z}\Q\to F\otimes_{\Z}\Q/\Z\to 0.\]
Apparently the vanishing of $R^i\varepsilon_{\mr{fl}*}F$ in Theorem \ref{1.2} is reduced to the vanishing of $H^i_{\mr{kfl}}(X,F)$ for any $X\in(\mr{fs}/S)$ such that the underlying scheme of $X$ is $\Spec R$ with $R$ a strictly henselian local ring. For the proof of Theorem \ref{1.1}, we first construct a canonical map
\[\Phi: j_{\mr{fl}!}R^i\varepsilon_{\mr{fl}*}j_{\mr{kfl}}^{-1}F\to R^i\varepsilon_{\mr{fl}*}F,\]
see (\ref{eqB.2}), then determine $R^i\varepsilon_{\mr{fl}*}j_{\mr{kfl}}^{-1}F$ (note that the order of $F$ is invertible on $U$), and finally prove that $\Phi$ is an isomorphism. The main computation tools are \v{C}ech cohomology, \v{C}ech-to-derived functor spectral sequence, and Leray spectral sequence.

In section \ref{sec3}, we apply the results from Section \ref{sec2} to make some computations on Kummer log flat cohomology when the base is a standard log trait or a Dedekind scheme endowed with the log structure associated to a finite set of closed points.
 
\subsection*{An application}
Theorem \ref{1.3} (1) can be used to prove that log abelian varieties with constant degeneration are sheaves for the Kummer log flat topology. According to \cite[Def. 3.3]{k-k-n2}, a log abelian variety with constant degeneration $A$ over an fs log scheme $S$ is a sheaf on $(\mr{fs}/S)_{\mr{\acute{e}t}}$ which is isomorphic to the quotient sheaf $G_{\mr{log}}^{(Y)}/Y$ for a pointwise polarizable log 1-motive $[Y\xrightarrow{u}G_{\mr{log}}]$ over $S$. One can consider $G_{\mr{log}}^{(Y)}$ as a uniformization of $A$, and $Y$ as the corresponding periods lattice. It is not surprising that $Y$ is important for understanding $A$. Indeed we have $R^1\varepsilon_{\mr{fl}*}Y=0$ by Theorem \ref{1.3} (1), and this vanishing is used in \cite[Thm. 2.1 (1)]{zha1} to prove that $A$ is a sheaf for the Kummer log flat topology.

In fact Theorem \ref{1.3} (1) is just \cite[Lem. 2.4]{zha1}. However the proof loc. cit. makes use of fpqc descent of schemes which probably does not always hold and deserves a precise reference (see \cite[\href{https://stacks.math.columbia.edu/tag/0APK}{Lemma 0APK}]{stacks-project} for the situation loc. cit.). The original motivation of this article is to present a new proof to \cite[Lem. 2.4]{zha1}.

\section{The higher direct images}\label{sec2}
We make a few lemmas first.
\begin{lem}\label{2.1}
Let $X$ be an fs log scheme. We make the following assumption on $X$, as well as some constructions associated to $X$.
\begin{itemize}
\item[$\bigstar$] The underlying scheme is $\Spec R$ with $R$ a strictly henselian local ring. Let $x$ denote the closed point, $k$ the residue field of $R$, $p$ the characteristic of $k$, and $P_X\to M_X$ a chart of the log structure of $X$ with $P$ an fs monoid, such that $P\xrightarrow{\cong}M_{X,x}/\mc{O}_{X,x}^{\times}$. Let $P^{1/n}$ denote the monoid $P$ regarded as a monoid above $P$ via the homomorphism $P\xrightarrow{n} P$. Let 
\[X_n:=X\times_{\Spec\Z[P]}\Spec\Z[P^{1/n}]\]
endowed with the canonical log structure associated to $P^{1/n}$ and let $H_n$ denote the group scheme $\Spec \Z[(P^{1/n})^{\mr{gp}}/P^{\mr{gp}}]$ over $\Spec\Z$, then $X_n$ is a Kummer log flat cover of $X$ such that 
\[X_n\times_XX_n\cong X_n\times_{\Spec\Z}H_n,\]
and it is even a Kummer log \'etale cover in case $(p,n)=1$.
\end{itemize}
Let $F$ be a constant group scheme over $\Spec R$ associated to an abelian group, and we write $n=m\cdot p^t$ with $(m,p)=1$. Then we have an isomorphism
\[\check{H}_{\mr{kfl}}^i(X_n/X,F)\cong H^i(H_m(X),F),\]
where the first term is the $i$-th \v{C}ech cohomology of $F$ with respect to the cover $X_n/X$, and the second term is the abstract group cohomology of the abstract group $H_m(X)$ with coefficients in $F$.

In particular, we have the following.
\begin{enumerate}
\item $\check{H}_{\mr{kfl}}^i(X_n/X,\Q)=0$ for $i>0$.
\item $\check{H}_{\mr{kfl}}^1(X_n/X,\Z)=0$ and $\check{H}_{\mr{kfl}}^i(X_n/X,\Z)$ is torsion and $p$-torsion-free for $i>0$.
\item The canonical map 
\[\check{H}_{\mr{kfl}}^i(X_m/X,\Z)\xrightarrow{\cong}\check{H}_{\mr{kfl}}^i(X_n/X,\Z)\]
is an isomorphism for $i>0$.
\item Suppose that $F$ is killed by a power of $p$. Then
\[\varinjlim_n\check{H}_{\mr{kfl}}^i(X_n/X,F)=0.\]
\item Suppose that $F$ is torsion and $p$-torsion-free. Then
\[\varinjlim_n\check{H}_{\mr{kfl}}^i(X_n/X,F)\cong\varinjlim_n H^i(H_m(X),F)\cong F(-i)\otimes_{\Z}\bigwedge^iP^{\mr{gp}}\]
for $i>0$. 
\end{enumerate}
\end{lem}
\begin{proof}
Since $R$ is a strictly henselian local ring, $X\times_{\Spec\Z}H_m^r$ is the constant group scheme over $X$ associated to the abstract group $H_m(X)^r$ and $X\times_{\Spec\Z}H_{p^t}^r$ is a connected group scheme over $X$, therefore we have
\begin{align*}
\Gamma(\underbrace{X_n\times_X\cdots\times_XX_n}_{\text{$r+1$ times}},F)=&\Gamma(X_n\times_{\Spec\Z}H_n^r,F)  \\
=&\Gamma((X\times_{\Spec\Z}H_{p^t}^r)\times_X(X_n\times_{\Spec\Z}H_m^r),F)\\
=&\prod_{h\in H_m(X)^r}F  \\
=&\mr{Map}(H_m(X)^r,F).
\end{align*}
We consider the \v{C}ech complex 
\begin{equation}\label{eq2.1}
\Gamma(X_n,F)\xrightarrow{d_0}\Gamma(X_n\times_XX_n,F)\xrightarrow{d_1}\Gamma(X_n\times_XX_n\times_XX_n,F)\xrightarrow{d_2}\cdots
\end{equation}
for $F$ with respect to the cover $X_n/X$. Let $\Gamma_n:=(P^{\frac{1}{n}})^{\mr{gp}}/P^{\mr{gp}}$. By \cite[Chap. III, Example 2.6]{mil1}, the \v{C}ech nerve of the Kummer log flat cover $X_n/X$ can be identified with the sequence
$$\xymatrix{
X_n & X_n\times H_n\ar@<0.5ex>[l]^-{d_{1,0}}\ar@<-0.5ex>[l]_-{d_{1,1}} &X_n\times H_n^2\ar@<1ex>[l]^-{d_{2,0}}\ar@<0ex>[l]\ar@<-1ex>[l]_-{d_{2,2}} &X_n\times H_n^3\ar@<1.5ex>[l]^-{d_{3,0}}\ar@<0.5ex>[l]\ar@<-0.5ex>[l]\ar@<-1.5ex>[l]_-{d_{3,3}}\cdots ,
}$$
where the map $d_{r,i}$ on the ring level is given by the $R$-linear ring homomorphism
\begin{align*}
R\otimes_{\Z[P]}\Z[P^{\frac{1}{n}}\oplus\Gamma_n^{r-1}]&\rightarrow R\otimes_{\Z[P]}\Z[P^{\frac{1}{n}}\oplus\Gamma_n^{r}]
  \\
(a,\bar{a}_1,\cdots,\bar{a}_{r-1})&\mapsto
\begin{cases} (a,\bar{a},\bar{a}_1,\cdots,\bar{a}_{r-1}),&\text{if $i=0$;}  \\
(a,\bar{a}_1,\cdots,\bar{a}_{i},\bar{a}_i,\cdots,\bar{a}_{r-1}),&\text{if $0<i<r$;}  \\
 (a,\bar{a}_1,\cdots,\bar{a}_{r-1},0),&\text{if $i=r$.}
\end{cases}
\end{align*}
for any $(a,\bar{a}_1,\cdots,\bar{a}_{r-1})\in P^{\frac{1}{n}}\oplus\Gamma_n^{r-1}$. If $m=1$, i.e. $n=p^t$, we have 
$$\Gamma(X_{p^t}\times_{\Spec\Z}H_{p^t}^r,F)=F.$$
The map 
$$d_{r,i}^*:\Gamma(X_{p^t}\times_{\Spec\Z}H_{p^t}^{r-1},F)\rightarrow \Gamma(X_{p^t}\times_{\Spec\Z}H_{p^t}^r,F)$$
is clearly just the identity map $\mr{Id}:F\rightarrow F$. In general, the map 
\[d_{r,i}^*:\Gamma(X_{n}\times_{\Spec\Z}H_{n}^{r-1},F)\rightarrow \Gamma(X_{n}\times_{\Spec\Z}H_{n}^r,F)\]
can be identified with the map
\[ \mr{Map}(H_m(X)^{r-1},F)\xrightarrow{\partial_{r,i}} \mr{Map}(H_m(X)^r,F)\] which maps $f\in \mr{Map}(H_m(X)^{r-1},F)$ to
\[\partial_{r,i}(f):(h_1,\cdots,h_r)\mapsto 
\begin{cases}
f(h_2,\cdots,\cdots,h_r), &\text{if $i=0$;} \\
f(h_1,\cdots,h_i+h_{i+1},\cdots,h_r), &\text{if $0<i<r$;} \\
f(h_1,\cdots\cdots,h_{r-1}), &\text{if $i=r$.} 
\end{cases}\]
Therefore the complex (\ref{eq2.1}) can be identified with the standard complex
\[F\rightarrow \mr{Map}(H_m(X),F)\rightarrow \mr{Map}(H_m(X)^2,F)\rightarrow\cdots\]
which computes the abstract group cohomology of the trivial $H_m(X)$-module $F$. It follows that $\check{H}_{\mr{kfl}}^i(X_n/X,F)\cong H^i(H_m(X),F)$.

In particular, we have part (1), part (2), and part (4) by the corresponding results for the abstract group cohomology $H^i(H_m(X),F)$. Part (3) is also clear by the construction of the isomorphism $\check{H}_{\mr{kfl}}^i(X_n/X,F)\cong H^i(H_m(X),F)$. 

We are left with part (5). We have
\begin{align*}
\varinjlim_n\check{H}_{\mr{kfl}}^i(X_n/X,F)\cong&\varinjlim_{\text{$n=mp^t$ with $(p,m)=1$}} H^i(H_m(X),F)  \\
\overset{(3)}{=}&\varinjlim_{(p,m)=1} H^i(H_m(X),F) \\
=&H^i(\varprojlim_{(p,m)=1}H_m(X),F)  \\
=&H^i(\mr{Hom}(P^{\mr{gp}},\hat{\Z}'(1)),F). \\
\end{align*}
By the description of the profinite group cohomology of $(\hat{\Z}')^r$ from Lemma \ref{C.1}, we further have
\begin{align*}
\varinjlim_n\check{H}_{\mr{kfl}}^i(X_n/X,F)\cong&H^i(\mr{Hom}(P^{\mr{gp}},\hat{\Z}'(1)),F) \\
=&\varinjlim_{(n,p)=1}H^i(\mr{Hom}(P^{\mr{gp}},\hat{\Z}'(1)),F[n])  \\
=&\varinjlim_{(n,p)=1}F[n]\otimes_{\Z/n\Z}H^i(\mr{Hom}(P^{\mr{gp}},\hat{\Z}'(1)),\Z/n\Z)  \\
=&\varinjlim_{(n,p)=1}F[n]\otimes_{\Z/n\Z}\bigwedge^iH^1(\mr{Hom}(P^{\mr{gp}},\hat{\Z}'(1)),\Z/n\Z) \\
=&\varinjlim_{(n,p)=1}F[n]\otimes_{\Z/n\Z}\bigwedge^i\mr{Hom}(\mr{Hom}(P^{\mr{gp}},\hat{\Z}'(1)),\Z/n\Z)  \\
=&\varinjlim_{(n,p)=1}F[n]\otimes_{\Z/n\Z}(\bigwedge^iP^{\mr{gp}}\otimes_{\Z}\Z/n\Z(-i))  \\
=&F(-i)\otimes_{\Z}\bigwedge^iP^{\mr{gp}}
\end{align*}
for $i>0$. This finishes the proof of part (5).
\end{proof}

\begin{rmk}
In the proof of Lemma \ref{2.1} (5), one can use Lemma \ref{C.1} to get $H^i(\mr{Hom}(P^{\mr{gp}},\hat{\Z}'(1)),F)\cong F\otimes_{\Z}\bigwedge^iP^{\mr{gp}}$ directly. However we repeat the proof of Lemma \ref{C.1} here in order to keep track of the Tate twist.
\end{rmk}

\begin{lem}\label{2.2}
Let $X$ and $X_n$ be as in Lemma \ref{2.1}. Let $F$ be a sheaf on $(\mr{fs}/X)_{\mr{kfl}}$. Then the family $\mathscr{X}_{\N}:=\{X_n\rightarrow X\}_{n\geq1}$ of Kummer log flat covers of $X$ satisfies the condition (L3) from \cite[\S 2]{art3}, whence a \v{C}ech-to-derived functor spectral sequence
\begin{equation}\label{eq2.2}
E_2^{i,j}=\check{H}_{\mr{kfl}}^i(\mathscr{X}_{\N},\underline{H}_{\mr{kfl}}^j(F))\Rightarrow H^{i+j}_{\mr{kfl}}(X,F),
\end{equation}
where 
$$\check{H}_{\mr{kfl}}^i(\mathscr{X}_{\N},-):=\varinjlim_{n\geq1}\check{H}_{\mr{kfl}}^i(X_n/X,-).$$
\end{lem}
\begin{proof}
This follows from \cite[Chap. II, Sec. 3, (3.3)]{art3}.
\end{proof}

A nice thing about the \v{C}ech cohomology with respect to all coverings, is that the zero-th cohomology vanishes, see \cite[Prop. 3.1]{zha5}. This is often very useful for computations. This is of course not the case for the \v{C}ech cohomology with respect to an arbitrary family. The following lemma is about the vanishing of the zero-th \v{C}ech cohomology for the covering family $\mathscr{X}_\N$ from Lemma \ref{2.2} under a suitable assumption.

\begin{lem}\label{2.3}
Let $X$ and $X_n$ be as in Lemma \ref{2.1}. Let $F$ be a sheaf on $(\mr{fs}/X)_{\mr{kfl}}$, and $N$ a non-negative integer. For any $U\in(\mr{fs}/X)$, let $(\mr{st}/U)$ be the full subcategory of $(\mr{fs}/U)$ consisting of strict objects over $U$. Assume that for any $U\in(\mr{fs}/X)$ and any $0\leq i\leq N$, the restriction of the flat sheaf $R^i\varepsilon_{\mr{fl}*}F$ to $(\mr{st}/U)$ is the inverse image of some sheaf on the small \'etale site of $U$, i.e. it lies in the image of the functor $a_U^{-1}$ from \cite[\href{https://stacks.math.columbia.edu/tag/0DDU}{Lemma 0DDU}]{stacks-project}.
Then we have 
\[\varinjlim_{n\geq1}H_{\mr{kfl}}^i(X_n,F)=0\]
for any $0<i\leq N+1$, in particular $\varinjlim_{n\geq1}\check{H}_{\mr{kfl}}^0(X_n/X,\underline{H}_{\mr{kfl}}^i(F))=0$ for any $0<i\leq N+1$.
\end{lem}
\begin{proof}
Clearly we only need to prove the vanishing for $i=N+1$. It suffices to show that for any $\gamma\in H_{\mr{kfl}}^{N+1}(X_r,F)$ there exists $s$ such that $\gamma$ goes to zero in $H_{\mr{kfl}}^{N+1}(X_{rs},F)$.  

By \cite[Prop. 3.1]{zha5}, we can find a Kummer log flat cover $T\rightarrow X_r$ such that $\gamma$ dies in $H_{\mr{kfl}}^{N+1}(T,F)$. By \cite[Cor. 2.16]{niz1}, we may assume that for some $s$, we have a factorization $T\rightarrow X_{rs}\rightarrow X_r$ with $T\rightarrow X_{rs}$ a classical flat cover. It follows that the pull-back $\gamma_s$ of $\gamma$ to $X_{rs}$ is trivialized by a classical flat cover, i.e. 
$$\gamma_s\in\mr{ker}(H_{\mr{kfl}}^{N+1}(X_{rs},F)\rightarrow H_{\mr{fl}}^0(X_{rs},R^{N+1}\varepsilon_{\mr{fl}*}F)).$$
By our assumption on the flat sheaf $R^i\varepsilon_{\mr{fl}*}F$ for any $0\leq i\leq N$, we get 
\[H^a_{\mr{fl}}(X_{rs},R^i\varepsilon_{\mr{fl}*}F)=H^a_{\mr{\acute{e}t}}(X_{rs},R^i\varepsilon_{\mr{fl}*}F)=0\]
for any $a>0$ and any $0\leq i\leq N$ by \cite[\href{https://stacks.math.columbia.edu/tag/0DDU}{Lemma 0DDU}]{stacks-project}. Therefore the Leray spectral sequence 
\[E_2^{u,v}=H^u_{\mr{fl}}(X_{rs},R^v\varepsilon_{\mr{fl}*}F)\Rightarrow H^{v+v}_{\mr{kfl}}(X_{rs},F),\]
implies that $\gamma_s=0$. This finishes the proof.
\end{proof}

In the setting-up of Lemma \ref{2.3}, the following lemma shows that $H^a_{\mr{kfl}}(X_{n,d},F)$ is determined only by $R^a\varepsilon_{\mr{fl}*}F$ for $0<a<N+1$, where $X_{n,d}$ denotes the fiber product of $d+1$ copies of $X_n$ over $X$. Apparently this is useful for computations on the spectral sequence (\ref{eq2.2}).

\begin{lem}\label{2.4}
Let the setting-up be as in Lemma \ref{2.3}, and let $X_{n,d}$ denote the fiber product of $d+1$ copies of $X_n$ over $X$. Then we have
\[H^a_{\mr{kfl}}(X_{n,d},F)=\Gamma(X_{n,d},R^a\varepsilon_{\mr{fl}*}F)\]
for $0<a<N+1$.
\end{lem}
\begin{proof}
We have $X_{n,d}=H_n^d\times_{\Spec\Z}X_n=(H_n)_{X}^d\times_XX_n$, where $H_n:=\Spec\Z[(P^{1/n})^{\mr{gp}}/P^{\mr{gp}}]$ and $(H_n)_X:=H_n\times_{\Spec\Z}X$. Note that $(H_n)_X$ is the constant group scheme associated to the abstract group $H_n(X)$ over $X$ if $(n,p)=1$. We write $n=p^r\cdot n'$ with $(p,n')=1$, then we have 
\[X_{n,d}=(H_{n'})_{X}^d\times_X((H_{p^r})_{X}^d\times_XX_n).\]
Thus 
\begin{align*}
H^u_{\mr{fl}}(X_{n,d},R^{v}\varepsilon_{\mr{fl}*}F)&=H^u_{\mr{\acute{e}t}}(X_{n,d},R^{v}\varepsilon_{\mr{fl}*}F) \\
&=\prod_{x\in H_{n'}(X)^d}H^u_{\mr{\acute{e}t}}((H_{p^r})_{X}^d\times_XX_n,R^{v}\varepsilon_{\mr{fl}*}F)  \\
&=0
\end{align*}
for $1\leq u$ and $0\leq v\leq N$.
Therefore the result follows from the Leray spectral sequence $E_2^{u,v}=H^u_{\mr{fl}}(X_{n,d},R^v\varepsilon_{\mr{fl}*}F)\Rightarrow H^{v+v}_{\mr{kfl}}(X_{n,d},F)$.
\end{proof}

\subsection{Case of finite abelian $l$-groups}
\begin{lem}\label{2.5}
Let $l$ be a prime number, and $S$ a locally noetherian fs log scheme on which $l$ is invertible. Let $F$ be a finite \'etale group scheme over the underlying scheme of $S$, we endow it with the induced log structure from $S$. Assume that $F$ is killed by $l^r$ for some positive integer $r$, then we have 
\[R^i\varepsilon_{\mr{fl}*}F\cong F(-i)\otimes_{\Z}\bigwedge^i(\Gml/\Gm)_{S_{\mr{fl}}}\]
for $i>0$.
\end{lem}
\begin{proof}
The proof is analogous to that of \cite[Thm. 2.4]{k-n1}. By Kato's theorem, see \cite[Thm. 4.1]{kat2}, we have an isomorphism
\[R^1\varepsilon_{\mr{fl}*}F\cong\mc{H}om(\Z/l^r\Z(1),F)\otimes_{\Z}(\Gml/\Gm)_{S_{\mr{fl}}}=F(-1)\otimes_{\Z}(\Gml/\Gm)_{S_{\mr{fl}}}.\]
Using cup-product, this isomorphism induces a homomorphism
\[\varphi_i:F(-i)\otimes_{\Z}\bigwedge^i(\Gml/\Gm)_{S_{\mr{fl}}}\to F\otimes_{\Z}R^i\varepsilon_{\mr{fl}*}\Z/l^r\Z\to R^i\varepsilon_{\mr{fl}*}F.\]

To finish the proof, it suffices to prove that $\varphi_i$ is an isomorphism for $i\geq1$ which we will proceed by induction.

\textbf{Base case.} The case $i=1$ is just Kato's theorem.

\textbf{Inductive step.} Let $N$ be a positive integer, and we assume that $\varphi_i$ is an isomorphism for any $1\leq i\leq N$. We are going to show that $\varphi_{N+1}$ is also an isomorphism. The map $\varphi_{N+1}$ induces a map
\[\phi_{N+1}:\Gamma(X,F(-N-1))\otimes_\Z\bigwedge^{N+1}P^\mr{gp}\rightarrow H^{N+1}_{\mr{kfl}}(X,F)\]
for any $X\in(\mr{fs}/S)$ satisfying the assumption $\bigstar$ in Lemma \ref{2.1}. It suffices to show that $\phi_{N+1}$ is an isomorphism. 

In the spectral sequence (\ref{eq2.2})
\[E_2^{u,v}=\varinjlim_n\check{H}_{\mr{kfl}}^{u}(X_n/X,\underline{H}_{\mr{kfl}}^v(F))\Rightarrow H^{u+v}_{\mr{kfl}}(X,F),\]
the terms $E_2^{u,v}$ with $u+v\leq N+1$ and $v\neq0$ (see the picture below) vanish (see the picture below) by Lemma \ref{2.6} below.
\[\begin{tikzpicture}
\draw [<->, thick] (5,0) -- (0,0) -- (0,5);
\foreach \Point in {(0,2),(0,3),(0,4), (1,2), (2,2), (1,3) ,(0,1), (1,1), (2,1), (3,1)}{
    \node at \Point {$\times$};}
\foreach \Point in {(0,0), (1,0), (2,0), (3,0), (4,0)}{
    \node at \Point {$\bullet$};}    
\node [below] at (5,0) {$u$};
\node [below] at (4,0) {$N+1$};
\node [below] at (3,0) {$$};
\node [below] at (2,0) {$$};
\node [below] at (1,0) {$1$};
\node [left] at (0,5) {$v$};
\node [left] at (0,4) {$N+1$};
\node [left] at (0,1) {$1$};
\node [below left] at (0,0) {$0$};
\end{tikzpicture}\]
Therefore the canonical map $\varinjlim_n\check{H}_{\mr{kfl}}^{N+1}(X_n/X,F)\xrightarrow{\cong}H^{N+1}_{\mr{kfl}}(X,F)$ is an isomorphism. Thus we can identify $\phi_{N+1}$ as a map 
\[\Gamma(X,F(-N-1))\otimes_\Z\bigwedge^{N+1}P^\mr{gp}\rightarrow \varinjlim_n\check{H}_{\mr{kfl}}^{N+1}(X_n/X,F),\]
and still call the resulting map $\phi_{N+1}$. Since both $\phi_{N+1}$ and the identification of Lemma \ref{2.1} (5) are constructed via cup-product, they agree. In particular $\phi_{N+1}$ is an isomorphism. This finishes the induction step.
\end{proof}

\begin{lem}\label{2.6}
Let the setting-up be as in the inductive step of the proof of Lemma \ref{2.5}. Then the terms $E_2^{u,v}$ with $u+v\leq N+1$ and $v\neq0$ vanish in the spectral sequence (\ref{eq2.2}).
\end{lem}
\begin{proof}
We have $E_2^{0,v}=0$ for $0<v\leq N+1$ by Lemma \ref{2.3}. We are left with the case that $uv\neq0$ and $u+v\leq N+1$. 

We denote by $X_{n,d}$ the fiber product of $d+1$ copies of $X_n$ over $X$. We have $X_{n,d}=H_n^d\times_{\Spec\Z}X_n=(H_n)_{X}^d\times_XX_n$, where $H_n:=\Spec\Z[(P^{1/n})^{\mr{gp}}/P^{\mr{gp}}]$ and $(H_n)_X:=H_n\times_{\Spec\Z}X$. Note that $(H_n)_X$ is the constant group scheme associated to the abstract group $H_n(X)$ over $X$ if $(n,p)=1$.

We compute the \v{C}ech cohomology group $\check{H}_{\mr{kfl}}^u(X_n/X,\underline{H}_{\mr{kfl}}^v(F))$ for $0<v< N+1$. Let $n=p^r\cdot n'$ with $(p,n')=1$, then we have 
$$X_{n,d}=(H_{n'})_{X}^d\times_X((H_{p^r})_{X}^d\times_XX_n)$$
and 
\begin{align*}
H_{\mr{kfl}}^v(X_{n,d},F)&=\prod_{a\in H_{n'}(X)^d}H_{\mr{kfl}}^v((H_{p^r})_{X}^d\times_XX_n,F)  \\
&=\mr{Map}(H_{n'}(X)^d,H_{\mr{kfl}}^v((H_{p^r})_{X}^d\times_XX_{n},F)).
\end{align*}
By the inductive hypothesis from the proof of Lemma \ref{2.5}, the assumption of Lemma \ref{2.3} is satisfied. Thus by Lemma \ref{2.4}, we have
\begin{align*}
H_{\mr{kfl}}^v((H_{p^r})_{X}^d\times_XX_{n},F)=&\Gamma((H_{p^r})_{X}^d\times_XX_{n},R^v\varepsilon_{\mr{fl}*}F)  \\
=&\Gamma((H_{p^r})_{X}^d\times_XX_{n},F(-v)\otimes_{\Z}\bigwedge^v(\Gml/\Gm)_{S_{\mr{fl}}})  \\
=&\Gamma(X,F(-v))\otimes_{\Z}\bigwedge^v(P^{\frac{1}{n}})^{\mr{gp}}  \\
=&H_{\mr{kfl}}^v(X_{n},F).
\end{align*}
It follows that $H_{\mr{kfl}}^v(X_{n,d},F)=\mr{Map}(H_{n'}(X)^d,H_{\mr{kfl}}^v(X_{n},F))$. The \v{C}ech complex for $\underline{H}_{\mr{kfl}}^v(F)$ with respect to the cover $X_{n}/X$ can be identified with the standard complex that computes the group cohomology of $H_{\mr{kfl}}^v(X_{n},F)$ regarded as a trivial $H_{n'}(X)$-module, hence we get $\check{H}_{\mr{kfl}}^u(X_{n}/X,\underline{H}_{\mr{kfl}}^v(F))=H^u(H_{n'}(X),H_{\mr{kfl}}^v(X_{n},F))$. Therefore
\begin{align*}
\varinjlim_{n}\check{H}_{\mr{kfl}}^u(X_n/X,\underline{H}_{\mr{kfl}}^v(F))=&\varinjlim_{n=p^r\cdot n'} H^u(H_{n'}(X),H_{\mr{kfl}}^v(X_{n},F))  \\
=&H^u(\varprojlim_{n'}H_{n'}(X),\varinjlim_{n=p^r\cdot n'}H_{\mr{kfl}}^v(X_{n},F)),
\end{align*}
where the second identification follows from \cite[\S 2, Prop. 8]{ser1}. By Lemma \ref{2.3} we have $\varinjlim_{n}H_{\mr{kfl}}^v(X_{n},F)=0$, and thus $E_2^{u,v}=\varinjlim_{n}\check{H}_{\mr{kfl}}^u(X_n/X,\underline{H}_{\mr{kfl}}^v(F))=0$ for and $0<v< N+1$. This finishes the proof of the lemma. 
\end{proof}

\begin{lem}\label{2.7}
Let $S$ be a locally noetherian fs log scheme. Let $l$ be a prime number, $U$ the open locus on $S$ where $l$ is invertible, and $j:U\hookrightarrow S$ the corresponding strict open immersion. Let $F$ be a finite \'etale group scheme over the underlying scheme of $S$ such that it is killed by a power of $l$. We endow $F$ with the induced log structure from $S$. Assume that the sheaf $R^i\varepsilon_{\mr{fl}*}F$ is supported on $U$ for $i>0$ (which we will show later), then we have the following.
\begin{enumerate}
\item $R^i\varepsilon_{\mr{fl}*}F\cong j_{\mr{fl!}}((j_{\mr{kfl}}^{-1}F)(-i)\otimes_{\Z}\bigwedge^i(\Gml/\Gm)_{U_{\mr{fl}}})$.
\item Let $X\in(\mr{fs}/S)$ satisfying the assumption $\bigstar$ in Lemma \ref{2.1} with $p=l$, then $H^a_{\mr{fl}}(X,R^i\varepsilon_{\mr{fl}*}F)=0$.
\end{enumerate}
\end{lem}
\begin{proof}
By Proposition \ref{B.4}, we get $R^i\varepsilon_{\mr{fl}*}F\cong j_{\mr{fl!}}R^i\varepsilon_{\mr{fl}*}j_{\mr{kfl}}^{-1}F$. By Lemma \ref{2.5}, we have $R^i\varepsilon_{\mr{fl}*}j_{\mr{kfl}}^{-1}F\cong (j_{\mr{kfl}}^{-1}F)(-i)\otimes_{\Z}\bigwedge^i(\Gml/\Gm)_{U_{\mr{fl}}}$. Then part (1) follows.

By the description of from part (1), the restriction of $R^i\varepsilon_{\mr{fl}*}F$ to $(\mr{st}/X)_{\mr{fl}}$ is the inverse image of some sheaf on the small \'etale site of $X$. Thus $H^a_{\mr{fl}}(X,R^i\varepsilon_{\mr{fl}*}F)=H^a_{\mr{\acute{e}t}}(X,R^i\varepsilon_{\mr{fl}*}F)$ by \cite[\href{https://stacks.math.columbia.edu/tag/0DDU}{Lemma 0DDU}]{stacks-project}. Let $x$ denote the closed point of $X$. By Gabber's theorem, see \cite[\href{https://stacks.math.columbia.edu/tag/09ZI}{Theorem 09ZI}]{stacks-project}, we get
\[H^a_{\mr{\acute{e}t}}(X,R^i\varepsilon_{\mr{fl}*}F)=H^a_{\mr{\acute{e}t}}(x,R^i\varepsilon_{\mr{fl}*}F).\]
Since $R^i\varepsilon_{\mr{fl}*}F$ is supported on $U$ and $x\notin U$, part (2) follows.
\end{proof}

\begin{lem}\label{2.8}
Let $S$ be a locally noetherian fs log scheme. Let $l$ be a prime number, and $U$ the open locus on $S$ where $l$ is invertible. Let $F$ be a finite \'etale group scheme over the underlying scheme of $S$, and we endow it with the induced log structure from $S$. Assume that $F$ is killed by a power of $l$, then the sheaf $R^i\varepsilon_{\mr{fl}*}F$ is supported on $U$ for $i>0$.
\end{lem}
\begin{proof}
We use induction on $i$. 

\textbf{Base case.} First of all we consider the case $i=1$. It suffices to show that $H^1_{\mr{kfl}}(X,F)=0$ for any $X\in(\mr{fs}/S)$ satisfying the assumption $\bigstar$ in Lemma \ref{2.1} with $p=l$. 

The spectral sequence (\ref{eq2.2}) gives rise to an exact sequence
\[0\to \varinjlim_{n\geq1}\check{H}_{\mr{kfl}}^1(X_n/X,F)\to H_{\mr{kfl}}^1(X,F)\to \varinjlim_{n\geq1}\check{H}_{\mr{kfl}}^0(X_n/X,\underline{H}_{\mr{kfl}}^1(F)).\]
By Lemma \ref{2.1} (4), we have $\varinjlim_{n\geq1}\check{H}_{\mr{kfl}}^1(X_n/X,F)=0$. By Lemma \ref{2.3}, we have $\varinjlim_{n\geq1}\check{H}_{\mr{kfl}}^0(X_n/X,\underline{H}_{\mr{kfl}}^1(F))=0$. Hence $H_{\mr{kfl}}^1(X,F)=0$.

\textbf{Inductive step.} Let $N$ be a positive integer, and we assume that $R^i\varepsilon_{\mr{fl}*}F$ is supported on $U$ for any $1\leq i\leq N$. We are going to prove that the same is true for $i=N+1$. It suffices to show that $H^{N+1}_{\mr{kfl}}(X,F)=0$ for any $X\in(\mr{fs}/S)$ satisfying the assumption $\bigstar$ in Lemma \ref{2.1} with $p=l$. 

Let $\gamma\in H^{N+1}_{\mr{kfl}}(X,F)$. By \cite[Prop. 3.1]{zha5}, we can find a Kummer log flat cover $T\rightarrow X$ such that $\gamma$ dies in $H_{\mr{kfl}}^{N+1}(T,F)$. By \cite[Cor. 2.16]{niz1}, we may assume that for some $m$, we have a factorization $T\rightarrow X_m\rightarrow X$ with $T\rightarrow X_m$ a classical flat cover. It follows that the pull-back $\gamma_m$ of $\gamma$ to $X_m$ is trivialized by a classical flat cover, i.e. 
$$\gamma_m\in\mr{ker}(H_{\mr{kfl}}^{N+1}(X_m,F)\rightarrow H_{\mr{fl}}^0(X_m,R^{N+1}\varepsilon_{\mr{fl}*}F)).$$
Note that $X_m$ also satisfies the assumption $\bigstar$ in Lemma \ref{2.1} with $p=l$.
Since $R^i\varepsilon_{\mr{fl}*}F$ is supported on $U$ for $0<i<N+1$ by induction hypothesis, we have 
\[H^a_{\mr{fl}}(X_m,R^i\varepsilon_{\mr{fl}*}F)=0\]
for any $a\geq0$ and any $0<i<N+1$ by Lemma \ref{2.7} (2). We also have 
\[H^{N+1}_{\mr{fl}}(X_m,F)=H^{N+1}_{\mr{\acute{e}t}}(X_m,F)=0.\]
Thus the Leray spectral sequence 
$$E_2^{i,j}=H^i_{\mr{fl}}(X_m,R^j\varepsilon_{\mr{fl}*}F)\Rightarrow H^{i+j}_{\mr{kfl}}(X_m,F)$$
implies that $\gamma_m=0$.

It follows that 
\[\gamma\in \mr{Ker}(H^{N+1}_{\mr{kfl}}(X,F)\to \varinjlim_n\check{H}^0_{\mr{kfl}}(X_n/X,\underline{H}^{N+1}_{\mr{kfl}}(F))).\]
Therefore to show that $\gamma$ is itself zero, it suffices to show that 
\[\varinjlim_n\check{H}^i_{\mr{kfl}}(X_n/X,\underline{H}^{N+1-i}_{\mr{kfl}}(F))=0\]
for $0<i\leq N+1$. By Lemma \ref{2.1} (4), we are left with the case $0<i<N+1$. Further it suffices to show that $\check{H}^i_{\mr{kfl}}(X_n/X,\underline{H}^{N+1-i}_{\mr{kfl}}(F))=0$ for any $0<i<N+1$ and any $n$.

Let $n=l^rn'$ with $(l,n')=1$, then we have $H_{n'}\times_{\Spec\Z}X_n$ is a constant group over $X_n$ and
\[\underbrace{X_n\times_X\cdots\times_XX_n}_\text{$d+1$ folded}=X_n\times_{\Spec\Z} H_n^d=\bigcup_{a\in H_{n'}(X)}X_n\times_{\Spec\Z} H_{l^r}^d\]
with $X_n\times_{\Spec\Z} H_{l^r}^d$ satisfying the assumption $\bigstar$ in Lemma \ref{2.1} with $p=l$. By the same argument showing $\gamma_m=0$, one can also show that 
\[H^{N+1-i}_{\mr{kfl}}(X_n\times_{\Spec\Z} H_{l^r}^d,F)=0\]
for any $0<i<N+1$ and any $d\geq0$. Thus 
\[H^{N+1-i}_{\mr{kfl}}(\underbrace{X_n\times_X\cdots\times_XX_n}_\text{$d+1$ folded},F)=0\]
for any $0<i<N+1$ and any $d\geq0$. It follows that $\check{H}^i_{\mr{kfl}}(X_n/X,\underline{H}^{N+1-i}_{\mr{kfl}}(F))=0$ for any $0<i<N+1$ and any $n$. Therefore $\gamma=0$, and thus $H^{N+1}_{\mr{kfl}}(X,F)=0$. This shows that $R^{N+1}\varepsilon_{\mr{fl}*}F$ is supported on $U$.

\textbf{Conclusion.} The lemma is proven.
\end{proof}

\begin{thm}\label{2.9}
Let $S$ be a locally noetherian fs log scheme. Let $l$ be a prime number, $U$ the open locus on $S$ where $l$ is invertible, and $j:U\hookrightarrow S$ the corresponding strict open immersion. Let $F$ be a finite \'etale group scheme over the underlying scheme of $S$, and we endow it with the induced log structure from $S$. Assume that $F$ is killed by a power of $l$, then we have
\[R^i\varepsilon_{\mr{fl}*}F\cong j_{\mr{fl!}}((j_{\mr{kfl}}^{-1}F)(-i)\otimes_{\Z}\bigwedge^i(\Gml/\Gm)_{U_{\mr{fl}}}).\]
\end{thm}
\begin{proof}
This follows from Lemma \ref{2.7} and Lemma \ref{2.8}.
\end{proof}

\subsection{Case of rational vector spaces}
In this subsection, we consider the case that $F$ is \'etale locally isomorphic to a finite dimensional $\Q$-vector space.

\begin{thm}\label{2.10}
Let $S$ be an fs log scheme. Let $F$ be a group scheme over the underlying scheme of $S$ which is \'etale locally representable by a finite dimensional $\Q$-vector space. Then we have $R^i\varepsilon_{\mr{fl}*}F=0$ for $i\geq 1$.
\end{thm}
\begin{proof}
It suffices to consider the case $F=\Q$. We use induction on $i$. 

\textbf{Base case.} First of all we consider the case $i=1$. It suffices to show that $H^1_{\mr{kfl}}(X,\Q)=0$ for any $X\in(\mr{fs}/S)$ which satisfies the assumption $\bigstar$ in Lemma \ref{2.1}. The spectral sequence (\ref{eq2.2}) gives rise to an exact sequence
\[0\to \varinjlim_{n\geq1}\check{H}_{\mr{kfl}}^1(X_n/X,\Q)\to H_{\mr{kfl}}^1(X,\Q)\to \varinjlim_{n\geq1}\check{H}_{\mr{kfl}}^0(X_n/X,\underline{H}_{\mr{kfl}}^1(\Q)).\]
By Lemma \ref{2.1} (1), we have $\varinjlim_{n\geq1}\check{H}_{\mr{kfl}}^1(X_n/X,\Q)=0$. By Lemma \ref{2.3}, we have $\varinjlim_{n\geq1}\check{H}_{\mr{kfl}}^0(X_n/X,\underline{H}_{\mr{kfl}}^1(\Q))=0$. Hence $H_{\mr{kfl}}^1(X,\Q)=0$, and thus $R^1\varepsilon_{\mr{fl}*}\Q=0$.

\textbf{Inductive step.} Let $N$ be a positive integer, and assume that $R^i\varepsilon_{\mr{fl}*}\Q=0$ for any $1\leq i\leq N$. We are going to prove that $R^{N+1}\varepsilon_{\mr{fl}*}\Q=0$. It suffices to show that $H^{N+1}_{\mr{kfl}}(X,\Q)=0$ for any $X\in(\mr{fs}/S)$ which satisfies the assumption $\bigstar$ in Lemma \ref{2.1}. 

By the inductive hypothesis and Lemma \ref{2.3}, we have 
\[\varinjlim_{n\geq1}\check{H}^0_{\mr{kfl}}(X_n/X,\underline{H}^{N+1}_{\mr{kfl}}(\Q)))=0.\]
By the spectral sequence (\ref{eq2.2}), to show that $H^{N+1}_{\mr{kfl}}(X,\Q)=0$, it suffices to show that 
\[\varinjlim_{n\geq1}\check{H}^i_{\mr{kfl}}(X_n/X,\underline{H}^{N+1-i}_{\mr{kfl}}(\Q))=0\]
for $0<i\leq N+1$. By Lemma \ref{2.1} (1), we are left with the case $0<i<N+1$. Further it suffices to show that $\check{H}^i_{\mr{kfl}}(X_n/X,\underline{H}^{N+1-i}_{\mr{kfl}}(\Q))=0$ for any $0<i<N+1$ and any $n$. Since $R^t\varepsilon_{\mr{fl}*}\Q=0$ for $1\leq t\leq N$ and
\[H^r_{\mr{fl}}(\underbrace{X_n\times_X\cdots\times_XX_n}_\text{$d+1$ folded},\Q)=H^r_{\mr{fl}}(X_n\times_{\Spec\Z} H_n^d,\Q)=H^r_{\mr{\acute{e}t}}(X_n\times_{\Spec\Z} H_n^d,\Q)=0\]
for any $r>0$ and any $d\geq0$, the Leray spectral sequence 
$$E_2^{u,v}=H^u_{\mr{fl}}(\underbrace{X_n\times_X\cdots\times_XX_n}_\text{$d+1$ folded},R^v\varepsilon_{\mr{fl}*}\Q)\Rightarrow H^{u+v}_{\mr{kfl}}(\underbrace{X_n\times_X\cdots\times_XX_n}_\text{$d+1$ folded},\Q)$$
implies that $H^{N+1-i}_{\mr{kfl}}(\underbrace{X_n\times_X\cdots\times_XX_n}_\text{$d+1$ folded},\Q)=0$ for any $0<i<N+1$ and any $d\geq0$. It follows that $\check{H}^i_{\mr{kfl}}(X_n/X,\underline{H}^{N+1-i}_{\mr{kfl}}(\Q))=0$ for any $0<i<N+1$ and any $n$. Therefore $H^{N+1}_{\mr{kfl}}(X,\Q)=0$. This finishes the proof of $R^{N+1}\varepsilon_{\mr{fl}*}\Q=0$.

%Let $\gamma\in H^{N+1}_{\mr{kfl}}(X,\Q)$. By \cite[Prop. 3.1]{zha5}, we can find a Kummer log flat cover $T\rightarrow X$ such that $\gamma$ dies in $H_{\mr{kfl}}^{N+1}(T,\Q)$. By \cite[Cor. 2.16]{niz1}, we may assume that for some $m$, we have a factorization $T\rightarrow X_m\rightarrow X$ with $T\rightarrow X_m$ a classical flat cover. It follows that $\gamma$ regarded as a class on $X_m$ is trivialized by a classical flat cover, i.e. 
%$$\gamma\in\mr{ker}(H_{\mr{kfl}}^{N+1}(X_m,\Q)\rightarrow H_{\mr{fl}}^0(X_m,R^{N+1}\varepsilon_{\mr{fl}*}\Q)).$$
%Since $R^i\varepsilon_{\mr{fl}*}\Q=0$ for $0<i<N+1$ by induction hypothesis and 
%$$H^{N+1}_{\mr{fl}}(X_m,\Q)=H^{N+1}_{\mr{\acute{e}t}}(X_m,\Q)=0,$$
%the Leray spectral sequence 
%$$E_2^{i,j}=H^i_{\mr{fl}}(X_m,R^j\varepsilon_{\mr{fl}*}\Q)\Rightarrow H^{i+j}_{\mr{kfl}}(X_m,\Q)$$
%implies that $\gamma$ becomes zero in $H^{N+1}_{\mr{kfl}}(X_m,\Q)$. It follows that 
%\[\gamma\in \mr{Ker}(H^{N+1}_{\mr{kfl}}(X,\Q)\to \check{H}^0_{\mr{kfl}}(\mathscr{X}_{\N},\underline{H}^{N+1}_{\mr{kfl}}(\Q))).\]
%Therefore to show that $\gamma$ is itself zero, it suffices to show that 
%\[\check{H}^i_{\mr{kfl}}(\mathscr{X}_{\N},\underline{H}^{N+1-i}_{\mr{kfl}}(\Q))=0\]
%for $0<i\leq N+1$.

\textbf{Conclusion.} The theorem is proven.
\end{proof}

\subsection{Case of finite rank torsion-free abelian groups}
\begin{thm}\label{2.11}
Let $S$ be a locally noetherian fs log scheme, and $F$ a group scheme over the underlying scheme of $S$ which is \'etale locally isomorphic to a finite rank free abelian group. Then we have the following.
\begin{enumerate}
\item $R^1\varepsilon_{\mr{fl}*}F=0$.
\item Let $i>1$. For each prime number $l$, let $U_l$ be the locus on $S$ on which $l$ is invertible and ${}_lj:U_l\hookrightarrow S$ the corresponding strict open immersion. Then 
\begin{align*}
R^i\varepsilon_{\mr{fl}*}F&\cong\bigoplus_{l\text{ prime}}R^{i-1}\varepsilon_{\mr{fl}*}(F\otimes_\Z\Q_l/\Z_l) \\
&\cong \bigoplus_{l\text{ prime}}{}_lj_{\mr{fl}!}({}_lj^{-1}F\otimes_\Z\Q_l/\Z_l(-i+1)\otimes_{\Z}\bigwedge^{i-1}(\Gml/\Gm)_{(U_{l})_{\mr{fl}}}).
\end{align*}
\end{enumerate} 
\end{thm}
\begin{proof}
Applying the functor $\varepsilon_{\mr{fl}*}$ to the short exact sequence 
\[0\to F\to F\otimes_{\Z}\Q\to \bigoplus_{l\text{ prime}}F\otimes_{\Z}\Q_l/\Z_l\to 0\]
of sheaves of abelian groups on $(\mr{fs}/S)_{\mr{kfl}}$, we get a long exact sequence
\begin{align*}
0\to &F\to F\otimes_{\Z}\Q\to \bigoplus_{l\text{ prime}}F\otimes_{\Z}\Q_l/\Z_l\to R^1\varepsilon_{\mr{fl}*}F\to R^1\varepsilon_{\mr{fl}*}(F\otimes_{\Z}\Q)\to\cdots \\
\to &R^{i-1}\varepsilon_{\mr{fl}*}(F\otimes_{\Z}\Q)\to \bigoplus_{l\text{ prime}}R^{i-1}\varepsilon_{\mr{fl}*}(F\otimes_{\Z}\Q_l/\Z_l)\to R^i\varepsilon_{\mr{fl}*}F \\
\to &R^{i}\varepsilon_{\mr{fl}*}(F\otimes_{\Z}\Q)\to\cdots .
\end{align*}
By Theorem \ref{2.10}, we get an exact sequence
\[0\to F\to F\otimes_{\Z}\Q\to \bigoplus_{l\text{ prime}}F\otimes_{\Z}\Q_l/\Z_l\to R^1\varepsilon_{\mr{fl}*}F\to0\]
and 
\begin{equation}\label{eq2.3}
\bigoplus_{l\text{ prime}}R^{i-1}\varepsilon_{\mr{fl}*}(F\otimes_{\Z}\Q_l/\Z_l)\xrightarrow{\cong} R^i\varepsilon_{\mr{fl}*}F
\end{equation}
for $i>1$. 

Since the map $F\otimes_{\Z}\Q\to \bigoplus_{l\text{ prime}}F\otimes_{\Z}\Q_l/\Z_l$ remains surjective on $(\mr{fs}/S)_{\mr{fl}}$, we get $R^1\varepsilon_{\mr{fl}*}F=0$.

Part (2) follows from the isomorphism (\ref{eq2.3}) and Theorem \ref{2.9}.
\end{proof}

\section{Examples}\label{sec3}

\subsection{Some general results}
\begin{lem}\label{3.1}
Let $X$ be an fs log scheme whose underlying scheme is locally noetherian, and let $F$ be a group scheme over the underlying scheme of $X$ which is \'etale locally isomorphic to a finite rank free abelian group. We endow $F$ with the induced log structure from $X$. Then we have 
\begin{enumerate}
\item an isomorphism $H^1_{\mr{fl}}(X,F)\xrightarrow{\cong}H^1_{\mr{kfl}}(X,F)$,
\item and an exact sequence
\begin{equation}\label{eq3.1}
0\to H^2_{\mr{fl}}(X,F)\rightarrow H^2_{\mr{kfl}}(X,F)\to H^0_{\mr{fl}}(X,R^2\varepsilon_{\mr{fl}*}F)\to H^3_{\mr{fl}}(X,F)\to H^3_{\mr{kfl}}(X,F).
\end{equation}
\end{enumerate}
\end{lem}
\begin{proof}
The results follow from the Leray spectral sequence
\[E_2^{i,j}=H^i_{\mr{fl}}(X,R^j\varepsilon_{\mr{fl}*}F)\Rightarrow H^{i+j}_{\mr{kfl}}(X,F)\]
and the vanishing of $R^1\varepsilon_{\mr{fl}*}F$ from Theorem \ref{2.11} (1).
\end{proof}

\begin{prop}
Let $X$ be an fs log scheme whose underlying scheme is noetherian and normal. Then we have 
\[H^1_{\mr{kfl}}(X,\Z)=0.\]
\end{prop}
\begin{proof}
This follows from the isomorphism from Lemma \ref{3.1} (1) and the vanishing of $H^1_{\mr{fl}}(X,\Z)$ from \cite[Chapter 2, Prop. 2.4.2]{ct-s1}.
\end{proof}

\begin{lem}\label{3.3}
Let $X$, $F$ be as in Lemma \ref{3.1}, and we further assume that the ranks of the stalks of the \'etale sheaf $M_X^{\mr{gp}}/\mc{O}_X^{\times}$ are at most one. Then the exact sequence (\ref{eq3.1}) can be extended further as
\begin{align*}
0&\to H^2_{\mr{fl}}(X,F)\to H^2_{\mr{kfl}}(X,F)\to H^0_{\mr{fl}}(X,R^2\varepsilon_{\mr{fl}*}F)  \\
&\to H^3_{\mr{fl}}(X,F)\to H^3_{\mr{kfl}}(X,F) \to H^1_{\mr{fl}}(X,R^2\varepsilon_{\mr{fl}*}F)\to\cdots  \\
&\to H^i_{\mr{fl}}(X,F)\to H^i_{\mr{kfl}}(X,F)\to H^{i-2}_{\mr{fl}}(X,R^2\varepsilon_{\mr{fl}*}F)\to\cdots .
\end{align*}
\end{lem}
\begin{proof}
By Theorem \ref{2.11} (1), we have $R^1\varepsilon_{\mr{fl}*}F=0$. By Theorem \ref{2.11} (2) and our assumption, the restriction of $R^i\varepsilon_{\mr{fl}*}F$ to $(\mr{st}/X)$ vanishes for $i>2$, where $(\mr{st}/X)$ denotes the full subcategory of $(\mr{fs}/X)$ consisting of strict fs log schemes over $X$. Then the result is an easy exercise on the spectral sequence $E_2^{i,j}=H^i_{\mr{fl}}(X,R^j\varepsilon_{\mr{fl}*}F)\Rightarrow H^{i+j}_{\mr{kfl}}(X,F)$.
\end{proof}

\begin{lem}\label{3.4}
Let $X$ be as in Lemma \ref{3.3}. Let $F$ be a finite \'etale group scheme over the underlying scheme of $X$, and we endow it with the induced log structure from $X$. Then we have a long exact sequence
\begin{align*}
0&\to H^1_{\mr{fl}}(X,F)\to H^1_{\mr{kfl}}(X,F)\to H^0_{\mr{fl}}(X,R^1\varepsilon_{\mr{fl}*}F)  \\
&\to H^2_{\mr{fl}}(X,F)\to H^2_{\mr{kfl}}(X,F) \to H^1_{\mr{fl}}(X,R^1\varepsilon_{\mr{fl}*}F)\to\cdots  \\
&\to H^i_{\mr{fl}}(X,F)\to H^i_{\mr{kfl}}(X,F)\to H^{i-1}_{\mr{fl}}(X,R^1\varepsilon_{\mr{fl}*}F)\to\cdots .
\end{align*}
\end{lem}
\begin{proof}
By Theorem \ref{2.9} and our assumption, the restriction of $R^i\varepsilon_{\mr{fl}*}F$ to $(\mr{st}/X)$ vanishes for $i>1$. Then the result is an easy exercise on the spectral sequence $E_2^{i,j}=H^i_{\mr{fl}}(X,R^j\varepsilon_{\mr{fl}*}F)\Rightarrow H^{i+j}_{\mr{kfl}}(X,F)$.
\end{proof}

\subsection{Discrete valuation rings}
Throughout this subsection, let $R$ be a discrete valuation ring with fraction field $K$ and residue field $k$. Assume that $k$ is of positive characteristic $p$. Let $\pi$ be a uniformizer of $R$, and we endow $X=\Spec R$ with the log structure associated to the homomorphism $\N\rightarrow R,1\mapsto\pi$. Let $x$ be the closed point of $X$ and $i$ the closed immersion $x\hookrightarrow X$, and we endow $x$ with the induced log structure from $X$. Let $\eta$ be the generic point of $X$ and $j$ the open immersion $\eta\hookrightarrow X$.

\begin{ex}\label{ex3.1}
Let $F$ be a group scheme over the underlying scheme of $X$ which is \'etale locally isomorphic to a finite rank free abelian group. We have 
\[H^1_{\mr{fl}}(X,F)\xrightarrow{\cong} H^1_{\mr{kfl}}(X,F)\] 
by Lemma \ref{3.1} (1). Clearly $X$ satisfies the assumption on log structure from Lemma \ref{3.3}, thus we get a long exact sequence
\begin{align*}
0&\to H^2_{\mr{fl}}(X,F)\to H^2_{\mr{kfl}}(X,F)\to H^0_{\mr{fl}}(X,R^2\varepsilon_{\mr{fl}*}F)  \\
&\to H^3_{\mr{fl}}(X,F)\to H^3_{\mr{kfl}}(X,F) \to H^1_{\mr{fl}}(X,R^2\varepsilon_{\mr{fl}*}F)\to\cdots  \\
&\to H^i_{\mr{fl}}(X,F)\to H^i_{\mr{kfl}}(X,F)\to H^{i-2}_{\mr{fl}}(X,R^2\varepsilon_{\mr{fl}*}F)\to\cdots .
\end{align*}
Apparently  the log structure of $X$ is supported on the closed point $x$ and the restriction of $(\Gml/\Gm)_{X_{\mr{fl}}}$ to $(\mr{st}/X)$ is isomorphic to $i_*\Z$. Since the only non-invertible prime on $X$ is $p$, the restriction of $R^2\varepsilon_{\mr{fl}*}F$ to $(\mr{st}/X)$ is isomorphic to
\[i_*(F\otimes_\Z(\Q/\Z)'(-1)),\]
where $(\Q/\Z)'$ denotes the prime to $p$ part of $\Q/\Z$. Thus we have 
\begin{align*}
H^u_{\mr{fl}}(X,R^2\varepsilon_{\mr{fl}*}F)&=H^u_{\mr{fl}}(x,F\otimes_\Z(\Q/\Z)'(-1))  \\
&=H^u_{\mr{\acute{e}t}}(x,F\otimes_\Z(\Q/\Z)'(-1)) \\
&=H^u(\mr{Gal}(k^{\mr{s}}/k),F\otimes_\Z(\Q/\Z)'(-1)),
\end{align*}
where $k^{\mr{s}}$ is a separable closure of $k$ and the second identification follows from \cite[\href{https://stacks.math.columbia.edu/tag/0DDU}{Lemma 0DDU}]{stacks-project}. We also have 
\begin{align*}
H^u_{\mr{fl}}(X,F)\cong H^u_{\mr{\acute{e}t}}(X,F)\cong H^u_{\mr{\acute{e}t}}(x,F)\cong H^u(\mr{Gal}(k^{\mr{s}}/k),F),
\end{align*}
where the first identification also follows from \cite[\href{https://stacks.math.columbia.edu/tag/0DDU}{Lemma 0DDU}]{stacks-project}. Thus we have
\[H^1_{\mr{kfl}}(X,F)\cong H^1(\mr{Gal}(k^{\mr{s}}/k),F), \] 
and the above exact sequence can be rewritten as
\begin{equation}\label{eq3.2}
\begin{split}
0&\to H^2(\mr{Gal}(k^{\mr{s}}/k),F)\to H^2_{\mr{kfl}}(X,F)\to H^0(\mr{Gal}(k^{\mr{s}}/k),F\otimes_\Z(\Q/\Z)'(-1))  \\
&\to H^3(\mr{Gal}(k^{\mr{s}}/k),F)\to H^3_{\mr{kfl}}(X,F) \to H^1(\mr{Gal}(k^{\mr{s}}/k),F\otimes_\Z(\Q/\Z)'(-1))\to\cdots  \\
&\to H^i(\mr{Gal}(k^{\mr{s}}/k),F)\to H^i_{\mr{kfl}}(X,F)\to H^{i-2}(\mr{Gal}(k^{\mr{s}}/k),F\otimes_\Z(\Q/\Z)'(-1))\to\cdots
\end{split}
\end{equation}

Assume that the cohomological dimension of $k$ is $d$, then we have
\[H^u(\mr{Gal}(k^{\mr{s}}/k),F)=0\] for $u>d+1$
and
\[H^u(\mr{Gal}(k^{\mr{s}}/k),F\otimes_\Z(\Q/\Z)'(-1))=0\]
for $u>d$. Thus we have an exact sequence
\begin{align*}
0&\to H^2(\mr{Gal}(k^{\mr{s}}/k),F)\to H^2_{\mr{kfl}}(X,F)\to H^0(\mr{Gal}(k^{\mr{s}}/k),F\otimes_\Z(\Q/\Z)'(-1))  \\
&\to H^3(\mr{Gal}(k^{\mr{s}}/k),F)\to H^3_{\mr{kfl}}(X,F) \to H^1(\mr{Gal}(k^{\mr{s}}/k),F\otimes_\Z(\Q/\Z)'(-1))\to\cdots  \\
&\to H^{d+1}(\mr{Gal}(k^{\mr{s}}/k),F)\to H^{d+1}_{\mr{kfl}}(X,F)\to H^{d-1}(\mr{Gal}(k^{\mr{s}}/k),F\otimes_\Z(\Q/\Z)'(-1))  \\
&\to0,
\end{align*}
an isomorphism 
\[H^{d+2}_{\mr{kfl}}(X,F)\xrightarrow{\cong}H^d(\mr{Gal}(k^{\mr{s}}/k),F\otimes_\Z(\Q/\Z)'(-1)),\]
and
\[H^i_{\mr{kfl}}(X,F)=0\]
for $i\geq d+3$ by the exact sequence (\ref{eq3.2}).

Moreover we assume that $k$ is finite, then $\mr{Gal}(k^{\mr{s}}/k)\cong\hat{\Z}$ and $d=1$. Thus we have an exact sequence
\begin{align*}
0&\to H^2(\hat{\Z},F)\to H^2_{\mr{kfl}}(X,F)\to H^0(\hat{\Z},F\otimes_\Z(\Q/\Z)'(-1)) \to 0,
\end{align*}
an isomorphism 
\[H^{3}_{\mr{kfl}}(X,F)\cong H^1(\hat{\Z},F\otimes_\Z(\Q/\Z)'(-1)),\]
and
\[H^i_{\mr{kfl}}(X,F)=0\]
for $i\geq 4$. Now we claim that
\[H^0(\hat{\Z},F\otimes_\Z(\Q/\Z)'(-1))=H^1(\hat{\Z},F\otimes_\Z(\Q/\Z)'(-1))=0.\]
This is clear if $F=\Z$. In general, take a finite extension $k'$ of $k$ such that $F\times_X\Spec k'\cong\Z^r$, then the claim follows from the Hochschild-Serre spectral sequence and the case of $F=\Z$. If follows that
\[H^i_{\mr{kfl}}(X,F)=\begin{cases}H^i(\hat{\Z},F),&\text{if $i=0,1,2$};\\
0, &\text{if $i>2$}.
\end{cases}\]
In particular we have
\[H^i_{\mr{kfl}}(X,\Z)=\begin{cases}\Z,&\text{if $i=0$};\\
0, &\text{if $i=1$}; \\
\Q/\Z, &\text{if $i=2$}; \\
0, &\text{if $i>2$}.
\end{cases}\]
\end{ex}

\begin{ex}
Let $F$ be a finite \'etale group scheme over the underlying scheme of $X$. Then $F=\bigoplus_{l\text{ prime}}F(l)$, where $F(l)$ is the $l$-primary subgroup of $F$ and also finite \'etale. To compute the Kummer log flat cohomology of $F$, it suffices to compute that of $F(l)$ for each prime $l$. Without loss of generality, we assume that $F=F(l)$ for some prime $l$. Since $X$ satisfies the assumption on log structure from Lemma \ref{3.3}, we get a long exact sequence
\begin{align*}
0&\to H^1_{\mr{fl}}(X,F)\to H^1_{\mr{kfl}}(X,F)\to H^0_{\mr{fl}}(X,R^1\varepsilon_{\mr{fl}*}F)  \\
&\to H^2_{\mr{fl}}(X,F)\to H^2_{\mr{kfl}}(X,F) \to H^1_{\mr{fl}}(X,R^1\varepsilon_{\mr{fl}*}F)\to\cdots  \\
&\to H^i_{\mr{fl}}(X,F)\to H^i_{\mr{kfl}}(X,F)\to H^{i-1}_{\mr{fl}}(X,R^1\varepsilon_{\mr{fl}*}F)\to\cdots .
\end{align*}
by Lemma \ref{3.4}.

Now we proceed by considering the two cases $l=p$ and $l\neq p$.

Case (1): $l=p$. In this case, we have $R^1\varepsilon_{\mr{fl}*}F=0$. Thus 
\[H^i_{\mr{kfl}}(X,F)\cong H^i_{\mr{fl}}(X,F)\cong H^i_{\mr{\acute{e}t}}(X,F)=H^i(\mr{Gal}(k^{\mr{s}}/k),F)\]
for any $i\geq0$.

Case (2): $l\neq p$.  Since $l\neq p$, we have $R^1\varepsilon_{\mr{fl}*}F\cong i_*F(-1)$. Similar to Example \ref{ex3.1}, we have 
\begin{align*}
H^u_{\mr{fl}}(X,R^1\varepsilon_{\mr{fl}*}F)=H^u_{\mr{fl}}(x,F(-1)) =H^u_{\mr{\acute{e}t}}(x,F(-1)) =H^u(\mr{Gal}(k^{\mr{s}}/k),F(-1))
\end{align*}
and
\begin{align*}
H^u_{\mr{fl}}(X,F)\cong H^u_{\mr{\acute{e}t}}(X,F)\cong H^u_{\mr{\acute{e}t}}(x,F)\cong H^u(\mr{Gal}(k^{\mr{s}}/k),F).
\end{align*}
Thus the above exact sequence can be rewritten as
\begin{equation}\label{eq3.3}
\begin{split}
0&\to H^1(\mr{Gal}(k^{\mr{s}}/k),F)\to H^1_{\mr{kfl}}(X,F)\to H^0(\mr{Gal}(k^{\mr{s}}/k),F(-1)) \\
&\to H^2(\mr{Gal}(k^{\mr{s}}/k),F)\to H^2_{\mr{kfl}}(X,F) \to H^1(\mr{Gal}(k^{\mr{s}}/k),F(-1))\to\cdots  \\
&\to H^i(\mr{Gal}(k^{\mr{s}}/k),F)\to H^i_{\mr{kfl}}(X,F)\to H^{i-1}(\mr{Gal}(k^{\mr{s}}/k),F(-1))\to\cdots .
\end{split}
\end{equation}
Assume that the cohomological dimension of $k$ is $d$, then we have
\[H^u(\mr{Gal}(k^{\mr{s}}/k),F)=H^u(\mr{Gal}(k^{\mr{s}}/k),F(-1))=0\] for $u>d$. Thus we have an exact sequence
\begin{align*}
0&\to H^1(\mr{Gal}(k^{\mr{s}}/k),F)\to H^1_{\mr{kfl}}(X,F)\to H^0(\mr{Gal}(k^{\mr{s}}/k),F(-1))  \\
&\to H^2(\mr{Gal}(k^{\mr{s}}/k),F)\to H^2_{\mr{kfl}}(X,F) \to H^1(\mr{Gal}(k^{\mr{s}}/k),F(-1))\to\cdots  \\
&\to H^{d}(\mr{Gal}(k^{\mr{s}}/k),F)\to H^{d}_{\mr{kfl}}(X,F)\to H^{d-1}(\mr{Gal}(k^{\mr{s}}/k),F(-1)) \to0,
\end{align*}
an isomorphism 
\[H^{d+1}_{\mr{kfl}}(X,F)\cong H^d(\mr{Gal}(k^{\mr{s}}/k),F(-1)),\]
and
\[H^i_{\mr{kfl}}(X,F)=0\]
for $i\geq d+2$ by the exact sequence (\ref{eq3.3}).

Moreover we assume that $k$ is finite, then $\mr{Gal}(k^{\mr{s}}/k)\cong\hat{\Z}$ and $d=1$. Thus we have an exact sequence
\[0\to H^1(\hat{\Z},F)\to H^1_{\mr{kfl}}(X,F)\to H^0(\hat{\Z},F(-1))\to0,\]
an isomorphism
\[H^{2}_{\mr{kfl}}(X,F)\cong H^1(\hat{\Z},F(-1)),\]
and
\[H^i_{\mr{kfl}}(X,F)=0\]
for $i\geq 3$. In particular, we have 
\[H^0(\hat{\Z},\Z/l^r\Z(-1))=H^1(\hat{\Z},\Z/l^r\Z(-1))=0,\]
and thus
\[H^i_{\mr{kfl}}(X,\Z/l^r\Z)=\begin{cases}\Z/l^r\Z,&\text{if $i=0$};\\
\Z/l^r\Z, &\text{if $i=1$}; \\
0, &\text{if $i>1$}.
\end{cases}\]
\end{ex}

\subsection{Global Dedekind domains}
Throughout this subsection, let $K$ be a global field. When $K$ is a number field, $X$ denotes the spectrum of the ring of integers in $K$, and when $K$ is a function field, $k$ denotes the field of constants of $K$ and $X$ denotes the unique connected smooth projective curve over $k$ having $K$ as its function field. Let $S$ be a finite set of closed points of $X$, $U:=X-S$, $j:U\hookrightarrow X$, and $i_x:x\hookrightarrow X$ for each closed point $x\in X$. We endow $X$ with the log structure $j_{*}\mc{O}^{\times}_U\cap\mc{O}_X\rightarrow \mc{O}_X$. 

Apparently $X$ satisfies the assumption on log structure from Lemma \ref{3.3}.

\begin{ex}\label{ex3.3}
Let $F$ be a group scheme over the underlying scheme of $X$ which is \'etale locally isomorphic to a finite rank free abelian group. We have 
\[H^1_{\mr{fl}}(X,F)\xrightarrow{\cong} H^1_{\mr{kfl}}(X,F)\] 
by Lemma \ref{3.1} (1). Since $X$ satisfies the assumption on log structure from Lemma \ref{3.3}, we have a long exact sequence as in Lemma \ref{3.3}. 

The sheaf $M_X/\mc{O}_X^\times$ is supported on the closed subset $S$, and thus the restriction of $(\Gml/\Gm)_{X_{\mr{fl}}}$ to $(\mr{st}/X)$ is isomorphic to $\bigoplus_{x\in S} i_{x*}\Z$. For any prime number $l$, let 
\[S_l:=\{x\in S\mid \text{the characteristic of the residue field of $x$ is not $l$}\}.\]
Thus by Theorem \ref{2.11},  we have
\[R^2\varepsilon_{\mr{fl}*}F\cong \bigoplus_{l\text{ prime}} R^1\varepsilon_{\mr{fl}*}(F\otimes_{\Z}\Q_l/\Z_l)\]
and
\[R^1\varepsilon_{\mr{fl}*}(F\otimes_{\Z}\Q_l/\Z_l)\cong\bigoplus_{x\in S_l}i_{x*}(F\otimes_{\Z}\Q_l/\Z_l(-1)).\]
Therefore we get
\begin{align*}
H^u_{\mr{fl}}(X,R^2\varepsilon_{\mr{fl}*}F)&=\bigoplus_{l\text{ prime}}\bigoplus_{x\in S_l}H^u_{\mr{fl}}(X,i_{x*}(F\otimes_{\Z}\Q_l/\Z_l(-1))) \\
&=\bigoplus_{l\text{ prime}}\bigoplus_{x\in S_l}H^u_{\mr{fl}}(x,F\otimes_{\Z}\Q_l/\Z_l(-1)) \\
&=\bigoplus_{l\text{ prime}}\bigoplus_{x\in S_l}H^u_{\mr{\acute{e}t}}(x,F\otimes_{\Z}\Q_l/\Z_l(-1)) \\
&=\bigoplus_{l\text{ prime}}\bigoplus_{x\in S_l}H^u(\Gamma_x,F\otimes_{\Z}\Q_l/\Z_l(-1)),
\end{align*}
where $\Gamma_x:=\mr{Gal}(\kappa(x)^{\mr{s}}/\kappa(x))$. We also have $H^2_{\mr{fl}}(X,F)\cong H^2_{\mr{\acute{e}t}}(X,F)$ by \cite[\href{https://stacks.math.columbia.edu/tag/0DDU}{Lemma 0DDU}]{stacks-project}. So we can rewrite the exact sequence from Lemma \ref{3.3} as
\begin{align*}
0&\to H^2_{\mr{\acute{e}t}}(X,F)\to H^2_{\mr{kfl}}(X,F)\to \bigoplus_{l\text{ prime}}\bigoplus_{x\in S_l}H^0(\Gamma_x,F\otimes_{\Z}\Q_l/\Z_l(-1))  \\
&\to H^3_{\mr{\acute{e}t}}(X,F)\to H^3_{\mr{kfl}}(X,F) \to \bigoplus_{l\text{ prime}}\bigoplus_{x\in S_l}H^1(\Gamma_x,F\otimes_{\Z}\Q_l/\Z_l(-1))\to\cdots  \\
&\to H^i_{\mr{\acute{e}t}}(X,F)\to H^i_{\mr{kfl}}(X,F)\to \bigoplus_{l\text{ prime}}\bigoplus_{x\in S_l}H^{i-2}(\Gamma_x,F\otimes_{\Z}\Q_l/\Z_l(-1))\to\cdots .
\end{align*}

Now we assume that the residue fields of $X$ at its closed points are finite, then we have 
\[H^u(\Gamma_x,F\otimes_\Z(\Q_l/\Z_l)(-1))=0\]
for $u>1$ due to cohomological dimension reason. Moreover, by the same argument as in Example \ref{ex3.1}, we even have 
\[H^0(\Gamma_x,F\otimes_\Z(\Q_l/\Z_l)(-1))=H^1(\Gamma_x,F\otimes_\Z(\Q_l/\Z_l)(-1))=0.\] 
It follows that
\[H^i_{\mr{kfl}}(X,F)\cong H^i_{\mr{\acute{e}t}}(X,F)\]
for any $i\geq0$.
\end{ex}

\begin{ex}
Let $F$ be a finite \'etale group scheme over the underlying scheme of $X$. Then $F=\bigoplus_{l\text{ prime}}F(l)$, where $F(l)$ is the $l$-primary subgroup of $F$ and also finite \'etale. To compute the Kummer log flat cohomology of $F$, it suffices to compute that of $F(l)$ for each prime $l$. Without loss of generality, we assume that $F=F(l)$ for some prime $l$. Similar to the situation of Example \ref{ex3.3}, we have
\[R^1\varepsilon_{\mr{fl}*}F\cong\bigoplus_{x\in S_l}i_{x*}(F(-1)).\]
Therefore by Lemma \ref{3.4} and similar arguments as in Example \ref{ex3.3}, we have a long exact sequence
\begin{align*}
0&\to H^1_{\mr{\acute{e}t}}(X,F)\to H^1_{\mr{kfl}}(X,F)\to \bigoplus_{x\in S_l} H^0(\Gamma_x,F(-1))  \\
&\to H^2_{\mr{\acute{e}t}}(X,F)\to H^2_{\mr{kfl}}(X,F) \to \bigoplus_{x\in S_l} H^1(\Gamma_x,F(-1))\to\cdots  \\
&\to H^i_{\mr{\acute{e}t}}(X,F)\to H^i_{\mr{kfl}}(X,F)\to \bigoplus_{x\in S_l} H^{i-1}(\Gamma_x,F(-1))\to\cdots .
\end{align*}
Now we assume that the residue fields of $X$ at the closed points are finite, then we have 
\[H^u(\Gamma_x,F(-1))=0\]
for $u>1$ due to cohomological dimension reason. Therefore we have an exact sequence
\begin{align*}
0&\to H^1_{\mr{\acute{e}t}}(X,F)\to H^1_{\mr{kfl}}(X,F)\to \bigoplus_{x\in S_l} H^0(\Gamma_x,F(-1))  \\
&\to H^2_{\mr{\acute{e}t}}(X,F)\to H^2_{\mr{kfl}}(X,F) \to \bigoplus_{x\in S_l} H^1(\Gamma_x,F(-1)) \\
&\to H^3_{\mr{\acute{e}t}}(X,F)\to H^3_{\mr{kfl}}(X,F)\to 0,
\end{align*}
and isomorphisms
\[H^i_{\mr{\acute{e}t}}(X,F)\xrightarrow{\cong} H^i_{\mr{kfl}}(X,F)\]
for $i\geq 4$. 

For $F=\Z/l^r\Z$, we further have
\[H^0(\Gamma_x,\Z/l^r\Z(-1))=H^1(\Gamma_x,\Z/l^r\Z(-1)) =0\]
for $x\in S_l$. Therefore 
\[H^i_{\mr{\acute{e}t}}(X,\Z/l^r\Z)\xrightarrow{\cong} H^i_{\mr{kfl}}(X,\Z/l^r\Z)\]
for $i\geq 0$. 
\end{ex}

\appendix
\section{Sites and sheaves}\label{app1}
In this appendix, we collect some general results about sites from \cite[\href{https://stacks.math.columbia.edu/tag/00UZ}{Chapter 00UZ}]{stacks-project}.
\subsection{Sites}
\begin{defn}
\cite[\href{https://stacks.math.columbia.edu/tag/00VH}{Definition 00VH}]{stacks-project} A \textbf{site} is given by a category $\mc{C}$ and a set $\mr{Cov}(\mc{C})=\bigcup_{U\in \mc{C}}\mr{Cov}(U)$ with $\mr{Cov}(U)$ being a set of families of morphisms with target $U$, such that the following conditions hold.
\begin{enumerate}
\item If $V\to U$ is an isomorphism, then $\{V\to U\}\in\mr{Cov}(U)$.
\item If $\{ U_ i \to U\} _{i\in I} \in \text{Cov}(U)$ and for each $i$ we have $\{ V_{ij} \to U_ i\} _{j\in J_ i} \in \text{Cov}(U_i)$, then $\{ V_{ij} \to U\} _{i \in I, j\in J_ i} \in \text{Cov}(U)$.
\item If $\{ U_ i \to U\} _{i\in I}\in \text{Cov}(U)$ and $V\to U$ is a morphism of $\mc{C}$, then $U_i\times_UV$ exists for each $i\in I$ and  $\{ U_ i \times _ U V \to V \} _{i\in I} \in \text{Cov}(V)$.
\end{enumerate}
\end{defn}

\subsection{Continuous functors}
\begin{defn}\cite[\href{https://stacks.math.columbia.edu/tag/00WV}{Tag 00WV}]{stacks-project}
A functor $u:\mc{C}\to\mc{D}$ of sites is called \textbf{continuous}, if for every $\{V_i\to V\}_{i\in I}\in\mr{Cov}(\mc{C})$ we have the following
 \begin{enumerate}
  \item $\{u(V_i)\to u(V)\}\in\mr{Cov}(\mc{D})$, and
  \item for any morphism $T\to V$ in $\mc{C}$ the morphism 
   \[u(T\times_VV_i)\to u(T)\times_{u(V)}u(V_i)\]
  is an isomorphism.
 \end{enumerate}
\end{defn}

Recall that given a functor $u$ as above, and a presheaf of sets $\mc{F}$ on $\mc{D}$ we can define $u^p\mc{F}$ to be simply the presheaf $\mc{F}\circ u$, in other words 
\[u^p\mc{F}(V)=\mc{F}(u(V))\]
for every object $V$ of $\mc{C}$ (see \cite[\href{https://stacks.math.columbia.edu/tag/00VC}{Tag 00VC}]{stacks-project} for $u^p$ as well as for $u_p$).

Suppose that the functor $u:\mc{C}\to\mc{D}$ is continuous, then $\mc{F}\in Sh(\mc{D})\Rightarrow u^p\mc{F}\in Sh(\mc{D})$. We denote
\[u^s:Sh(\mc{D})\to Sh(\mc{C})\]
the functor $u^p$ restricted to the subcategory of sheaves of sets. Recall that $u^p$ admits a left adjoint $u_p$, see \cite[\href{https://stacks.math.columbia.edu/tag/00VE}{Tag 00VE}]{stacks-project}. This is also the case for $u^s$.

\begin{lem}\cite[\href{https://stacks.math.columbia.edu/tag/00WX}{Tag 00WX}]{stacks-project}
Let $u:\mc{C}\to\mc{D}$ be a continuous functor of sites. Then the functor 
\[u_s:Sh(\mc{C})\to Sh(\mc{D}),\quad \mc{G}\mapsto (u_p\mc{G})^\sharp\]
is left adjoint to $u^s$.
\end{lem}

\begin{defn}\cite[\href{https://stacks.math.columbia.edu/tag/00X1}{Tag 00X1}]{stacks-project}
A \textbf{morphism of sites} $f:\mc{D}\to\mc{C}$ is given by a continuous functor $u:\mc{C}\to\mc{D}$ such that $u_s$ is exact.
\end{defn}
Notice how the functor $u$ goes in the direct opposite the morphism $f$. If $f\leftrightarrow u$ is a morphism of sites, then we use the notation $f^{-1}=u_s$ and $f_*=u^s$. The functor $f^{-1}$ is called the \textbf{pullback functor}, and the functor $f_*$ is called the \textbf{pushforward functor}. As in topology we have the adjunction $(f^{-1},f_*)$. 

See \cite[\href{https://stacks.math.columbia.edu/tag/00X2}{Tag 00X2}]{stacks-project} (examples associated to maps between two topological spaces) and \cite[\href{https://stacks.math.columbia.edu/tag/0EWI}{Tag 0EWI}]{stacks-project} (examples for different topologies on the same space, comparison of topologies) for examples of morphisms of sites.

\subsection{Cocontinuous functors}
\begin{defn}
A functor $u:\mc{C}\to\mc{D}$ of sites is called \textbf{cocontinuous}, if for every $U\in\mc{C}$ and every $\{V_j\to u(U)\}_{j\in J}\in\mr{Cov}(\mc{D})$, there exists $\{U_i\to U\}_{i\in I}\in\mr{Cov}(\mc{C})$ such that $\{u(U_i)\to u(U)\}_{i\in I}$ refines $\{V_j\to u(U)\}_{j\in J}$.
\end{defn}
\textbf{Warning}: In general $\{u(U_i)\to u(U)\}_{i\in I}$ is not a covering of $\mc{D}$.

\begin{ex}
For an fs log scheme $S$, we denote by $(\mr{Sch}/S)$ the category of fs log schemes over $S$. Let $(\mr{Sch}/S)_{\mr{k\acute{e}t}}$ and $(\mr{Sch}/S)_{\mr{\acute{e}t}}$ be the \textbf{Kummer log \'etale site} and the \textbf{classical \'etale site} for $(\mr{Sch}/S)$ respectively, see \cite[\S 2.5]{ill1}. Let $(\mr{Sch}/S)_{\mr{kfl}}$ be the \textbf{Kummer log flat site} for $(\mr{Sch}/S)$, see \cite[Def. 2.3]{kat2}, and let $(\mr{Sch}/S)_{\mr{fl}}$ be the \textbf{classical flat site} for $(\mr{Sch}/S)$, which is an obvious analogue of $(\mr{Sch}/S)_{\mr{\acute{e}t}}$.

Now let $j:U\hookrightarrow X$ be a strict open immersion of fs log schemes, and let $\tau\in\{\text{k\'et, \'et, kfl, fl}\}$. Then the functor of sites  
\[u:(\mr{fs}/X)_{\tau}\to(\mr{fs}/U)_{\tau},Y\mapsto Y\times_XU\]
is continuous, and the functor of sites
\[v:(\mr{fs}/U)_{\tau}\to (\mr{fs}/X)_{\tau},V\mapsto V\]
is continuous and cocontinuous
\end{ex}

\section{The Kummer log flat site}
In this appendix, we focus on the Kummer log flat site. 
\subsection{Morphisms of sites associated to a strict open immersion}
\begin{flushleft}
\textbf{Warning:}
We are going to follow \cite[\href{https://stacks.math.columbia.edu/tag/00UZ}{Chapter 00UZ}]{stacks-project} and \cite[\href{https://stacks.math.columbia.edu/tag/03A4}{Chapter 03A4}]{stacks-project} to construct a homomorphism of sheaves of abelian groups. Morphisms of sites and topoi in \cite[\href{https://stacks.math.columbia.edu/tag/00UZ}{Chapter 00UZ}]{stacks-project} are often formulated for sheaves of sets, while the map to be constructed is for sheaves of abelian groups. Often in order to use results from \cite[\href{https://stacks.math.columbia.edu/tag/00UZ}{Chapter 00UZ}]{stacks-project}, we have to use \cite[\href{https://stacks.math.columbia.edu/tag/00YV}{Tag 00YV}]{stacks-project}. Below whenever we refer to a result about sheaves of sets from \cite[\href{https://stacks.math.columbia.edu/tag/00UZ}{Chapter 00UZ}]{stacks-project} for sheaves of abelian groups, we are referring to \cite[\href{https://stacks.math.columbia.edu/tag/00YV}{Tag 00YV}]{stacks-project} without mention at the same time.
\end{flushleft}

Let $X$ be an fs log scheme, and $j:U\hookrightarrow X$ a strict open immersion of fs log schemes. Let $\tau\in\{\text{k\'et, \'et, kfl, fl}\}$. The functor
\[u:(\mr{fs}/U)_{\tau}\to (\mr{fs}/X)_{\tau},(V\to U)\mapsto (V\to U\to X)\]
of sites is continuous and cocontinuous, hence it gives rise to a morphism of topoi
\[g=(g^{-1},g_*):\mr{Ab}((\mr{fs}/U)_{\tau})\to \mr{Ab}((\mr{fs}/X)_{\tau})\]
with $g^{-1}$ exact by \cite[\href{https://stacks.math.columbia.edu/tag/00XO}{Tag 00XO}]{stacks-project} and \cite[\href{https://stacks.math.columbia.edu/tag/00XL}{Tag 00XL}]{stacks-project}.
For any presheaf $\mc{F}\in \mr{PAb}((\mr{fs}/U)_{\tau})$, we define $g_{p!}\mc{F}\in \mr{PAb}((\mr{fs}/X)_{\tau})$ as the presheaf
\[Y\mapsto \varinjlim_{Y\to u(V)}\mc{F}(V)\]
with colimits over $(\mc{I}^v_Y)^{\mr{opp}}$ (see \cite[\href{https://stacks.math.columbia.edu/tag/053L}{Equation 053L}]{stacks-project} for this index category) taken in the category of abelian groups. 
For $\mc{F}\in \mr{Ab}((\mr{fs}/U)_{\tau})$, we set $g_!\mc{F}$ to be the sheafification of $g_{p!}\mc{F}$, see \cite[\href{https://stacks.math.columbia.edu/tag/04BF}{Tag 04BF}]{stacks-project}, and called it the \textbf{extension by zero of $\mc{F}$}.

\begin{lem}
Let $X$, $U$, $u$, $g=(g^{-1},g_*)$, and $g_!$ be as above. Then we have the following.
\begin{enumerate}
\item The functor $g_!$ is left adjoint to $g^{-1}$ and exact.
\item The functor $g^{-1}$ is exact and preserves injective objects.
\item For any $\mc{F}\in\mr{Ab}((\mr{fs}/U)_{\tau})$, the canonical maps 
\[\mc{F}\to g^{-1}g_!\mc{F} \text{  and  } g^{-1}g_*\mc{F}\to \mc{F}\]
are isomorphisms.
\end{enumerate}
\end{lem}
\begin{proof}
Part (1) follows from \cite[\href{https://stacks.math.columbia.edu/tag/04BG}{Tag 04BG}]{stacks-project} and \cite[\href{https://stacks.math.columbia.edu/tag/04BH}{Tag 04BH}]{stacks-project}.

Since $g^{-1}$ admits both a left adjoint $g_!$ and a right adjoint $g_*$, it is exact by \cite[\href{https://stacks.math.columbia.edu/tag/0039}{Tag 0039}]{stacks-project}. Since the left adjoint $g_!$ of $g^{-1}$ is exact, it preserves injective objects by \cite[\href{https://stacks.math.columbia.edu/tag/015Z}{Tag 015Z}]{stacks-project}. This finishes the proof of part (2).

Since $j$ is a strict open immersion, the functor $u$ is fully faithful. The functor $u$ is continuous and cocontinuous. Thus part (3) follows from \cite[\href{https://stacks.math.columbia.edu/tag/077I}{Lemma 077I}]{stacks-project}.
\end{proof}

The functor 
\[v:(\mr{fs}/X)_{\tau}\to (\mr{fs}/U)_{\tau},Y\mapsto Y\times_XU\]
of sites is continuous with $v_s$ exact, hence gives rise to a morphism
\[j_{\tau}:(\mr{fs}/U)_{\tau}\to (\mr{fs}/X)_{\tau}\]
of sites, and further a morphism
\[j_{\tau}=(j_{\tau}^{-1}=v_s,j_{\tau*}=v^s):\mr{Ab}((\mr{fs}/U)_{\tau})\to \mr{Ab}((\mr{fs}/X)_{\tau})\]
of topoi. The cocontinuous functor $u$ is left adjoint to the continuous functor $v$. By \cite[\href{https://stacks.math.columbia.edu/tag/00XY}{Tag 00XY}]{stacks-project}, the two morphisms $g$ and $j_{\tau}$ of topoi agree. We set $j_{\tau!}:=g_!$, and call it the \textbf{functor of extension by zero}. To sum up, we have the following lemma.

\begin{lem}\label{B.2}
Let $j:U\to X$ be a strict open immersion of fs log schemes. Then we have a morphism of sites
\[j_{\tau}:(\mr{fs}/U)_{\tau}\to (\mr{fs}/X)_{\tau},\]
and a sequence of functors
\[j_{\tau!},j_{\tau}^{-1},j_{\tau*}\]
where in each consecutive pair the first is exact and left adjoint to the second. Moreover we have the following.
\begin{enumerate}
\item The functor $j_{\tau !}$ is exact.
\item The functor $j_{\tau}^{-1}$ is exact and preserves injective objects.
\item For any $\mc{F}\in\mr{Ab}((\mr{fs}/U)_{\tau})$, the canonical maps 
\[\mc{F}\to j_{\tau}^{-1}j_{\tau !}\mc{F} \text{ and } j_{\tau}^{-1}j_{\tau *}\mc{F}\to \mc{F}\]
are isomorphisms.
\end{enumerate}
\end{lem}

\subsection{Comparison morphism from the Kummer log flat site to the classical flat site}
Let $j:U\to X$ be a strict open immersion of fs log schemes. We have canonical forgetful morphisms of sites
\[(\mr{fs}/X)_{\mr{kfl}}\to (\mr{fs}/X)_{\mr{fl}}\]
and
\[(\mr{fs}/U)_{\mr{kfl}}\to (\mr{fs}/U)_{\mr{fl}}.\]
By abuse of notation, we denote both of them by $\varepsilon_{\mr{fl}}$. It is easy to see that we have the following commutative diagram
\[\xymatrix{
(\mr{fs}/U)_{\mr{kfl}}\ar[r]^{j_{\mr{kfl}}}\ar[d]_{\varepsilon_{\mr{fl}}} &(\mr{fs}/X)_{\mr{kfl}}\ar[d]^{\varepsilon_{\mr{fl}}}  \\
(\mr{fs}/U)_{\mr{fl}}\ar[r]_{j_{\mr{fl}}} &(\mr{fs}/X)_{\mr{fl}}
}\]
of morphisms of sites.

For $F\in\mr{Ab}((\mr{fs}/X)_{\mr{kfl}})$, the adjunction $(j_{\mr{fl}!},j_{\mr{fl}}^{-1})$ gives a canonical map
\begin{equation}\label{eqB.1}
j_{\mr{fl}!} j_{\mr{fl}}^{-1}R^i\varepsilon_{\mr{fl}*}F\to R^i\varepsilon_{\mr{fl}*}F.
\end{equation}

\begin{lem}\label{B.3}
We have 
\[j_{\mr{fl}}^{-1}R^i\varepsilon_{\mr{fl}*}F=R^i\varepsilon_{\mr{fl}*}j_{\mr{kfl}}^{-1}F.\]
\end{lem}
\begin{proof}
Let $F\to I^\bullet$ be an injective resolution of $F$ on $(\mr{fs}/X)_{\mr{kfl}}$. Since the functor $j_{\mr{kfl}}^{-1}$ is exact and preserves injective objects by Lemma \ref{B.2} (2), we get an injective resolution $j_{\mr{kfl}}^{-1}F\to j_{\mr{kfl}}^{-1}I^\bullet$ of $j_{\mr{kfl}}^{-1}F$. Thus 
\[R^i\varepsilon_{\mr{fl}*}j_{\mr{kfl}}^{-1}F=H^i(\varepsilon_{\mr{fl}*}j_{\mr{kfl}}^{-1}I^\bullet)=H^i(j_{\mr{fl}}^{-1}\varepsilon_{\mr{fl}*}I^\bullet)=j_{\mr{fl}}^{-1}H^i(\varepsilon_{\mr{fl}*}I^\bullet)=j_{\mr{fl}}^{-1}R^i\varepsilon_{\mr{fl}*}F.\]
\end{proof}

By the identification from Lemma \ref{B.3}, we get a map
\begin{equation}\label{eqB.2}
\Phi: j_{\mr{fl}!}R^i\varepsilon_{\mr{fl}*}j_{\mr{kfl}}^{-1}F\to R^i\varepsilon_{\mr{fl}*}F.
\end{equation}

\begin{prop}\label{B.4}
Let $X$ be an fs log scheme and $U$ an open subscheme of the underlying scheme of $X$. We endow $U$ with the induced log structure, and let $j:U\hookrightarrow X$ be the corresponding strict open immersion. Let $F$ be a sheaf of abelian groups on $(\mr{fs}/X)_{\mr{kfl}}$ such that $R^i\varepsilon_{\mr{fl}*}F$ is supported over $U$, then the map
\[\Phi: j_{\mr{fl}!}R^i\varepsilon_{\mr{fl}*}j_{\mr{kfl}}^{-1}F\to R^i\varepsilon_{\mr{fl}*}F\]
is an isomorphism.
\end{prop}
\begin{proof}
Applying the functor $j_{\mr{fl}}^{-1}$ to the canonical map (\ref{eqB.1}), we get a map
\begin{equation}\label{eqB.3}
j_{\mr{fl}}^{-1}j_{\mr{fl}!} j_{\mr{fl}}^{-1}R^i\varepsilon_{\mr{fl}*}F\to j_{\mr{fl}}^{-1}R^i\varepsilon_{\mr{fl}*}F.
\end{equation}
Since $j_{\mr{fl}}^{-1}j_{\mr{fl}!}j_{\mr{fl}}^{-1}R^i\varepsilon_{\mr{fl}*}F$ is identified to $j_{\mr{fl}}^{-1}R^i\varepsilon_{\mr{fl}*}F$ by Lemma \ref{B.2} (3), the map (\ref{eqB.3}) is an isomorphism. Since $R^i\varepsilon_{\mr{fl}*}F$ is supported on $U$ and the extension by zero sheaf $j_{\mr{fl}!} j_{\mr{fl}}^{-1}R^i\varepsilon_{\mr{fl}*}F$ is clearly supported on $U$, the map (\ref{eqB.1}) is actually an isomorphism. Therefore $\Phi$ is also an isomorphism by construction.
\end{proof}

\section{A lemma on profinite group cohomology}\label{app3}
Let $p$ be a fixed prime number, and let $\hat{\Z}':=\varprojlim_{(p,m)=1}\Z/m\Z$. Let $G_r=(\hat{\Z}')^r$ and $M$ a torsion abelian group, we regard $M$ as a $G_r$-module with respect to the trivial action. In this appendix, we compute the profinite group cohomology $H^i(G_r,M)$. The result should be well-known to the experts, but we are not able to find a reference so present a computation here. 

According to \cite[A.2]{ols2}, the cohomology groups $H^i(G_r,M)$ are computed by the cohomology groups of the standard homogeneous cochain complex of $G_r$ with coefficients in $M$
\[R\Gamma(G_r,M):\mr{Map}_{G_r}^{\mr{cts}}(G_r,M)\to \mr{Map}_{G_r}^{\mr{cts}}(G_r^2,M)\to\cdots\to \mr{Map}_{G_r}^{\mr{cts}}(G_r^{i},M)\to \cdots,\]
where $\mr{Map}_{G_r}^{\mr{cts}}(G_r^i,M)$ denotes the set of equivariant continuous functions $\phi:G_r^{i}\to M$ (where $M$ is endowed with the discrete topology). Note that $\mr{Map}_{G_r}^{\mr{cts}}(G_r^i,M)$ is denoted as $\mr{Hom}_{G_r}^{\mr{cts}}(G_r^{[i-1]},M)$ in \cite[App. A]{ols2}.

First we consider the case that $M=\Z/n\Z$ with $(n,p)=1$. The complex $R\Gamma(G_r,\Z/n\Z)$ of abelian groups is also naturally a complex of modules over the ring $\Z/n\Z$. Since $\Z/n\Z$ as a module over itself is flat and $R\Gamma(G_r,\Z/n\Z)$ lies in $D^b(\Z/n\Z)$ (the bounded derived category of complexes of $\Z/n\Z$-modules), we have
\begin{equation}\label{eqC.1}
R\Gamma(G_r,\Z/n\Z)\otimes_{\Z/n\Z}^L R\Gamma(G_s,\Z/n\Z)\xrightarrow{\cong}R\Gamma(G_{r+s},\Z/n\Z)
\end{equation}
by K\"unneth formula, see \cite[Thm. A.6]{ols2}. We have \[H^i(G_1,\Z/n\Z)=\begin{cases}\Z/n\Z,&\text{ if $i=0,1$;} \\ 0,&\text{ if $i>0$,}\end{cases}\]
and $R\Gamma(G_1,\Z/n\Z)$ is isomorphic to $\Z/n\Z\xrightarrow{0}\Z/n\Z$ in $D^b(\Z/n\Z)$. By \cite[Cor. A.7]{ols2}, the natural map of graded $\Z/n\Z$-modules
\[H^*(G_1,\Z/n\Z)\otimes_{\Z/n\Z}H^*(G_1,\Z/n\Z)\to H^*(G_2,\Z/n\Z)\]
is an isomorphism, and thus $H^i(G_2,\Z/n\Z)$ are free $\Z/n\Z$-modules for all $i$. Applying \cite[Cor. A.7]{ols2} inductively, we have that the natural map of graded $\Z/n\Z$-modules
\[H^*(G_1,\Z/n\Z)^{\otimes r}\to H^*(G_r,\Z/n\Z)\]
is an isomorphism. It follows that the graded cohomology ring $H^*(G_r,\Z/n\Z)$ is isomorphic to the exterior algebra of the module 
\[H^1(G_r,\Z/n\Z)=\mr{Hom}(G_r,\Z/n\Z)\cong (\Z/n\Z)^r\]
over $\Z/n\Z$.

Now let $M$ be a torsion abelian group which is killed by $n$ with $(n,p)=1$. We regard $M$ as a module over $\Z/n\Z$. Clearly we have 
\[R\Gamma(G_r,M)\cong R\Gamma(G_r,\Z/n\Z)\otimes_{\Z/n\Z}M.\]
Since $H^i(G_r,\Z/n\Z)$ are free $\Z/n\Z$-modules for all $i$, we get
\[H^i(G_r,M)\cong H^i(G_r,\Z/n\Z)\otimes_{\Z/n\Z}M\]
canonically.

In general for a torsion abelian group $M$, we have 
\[M=M[p^\infty]\oplus M'=M[p^\infty]\oplus\varinjlim_{(n,p)=1}M[n],\]
where $M[p^\infty]$ (resp. $M'$, resp. $M[n]$) denotes the $p$-primary part (resp. prime to $p$ part, resp. $n$-torsion part) of $M$. Thus
\begin{equation}\label{eqC.2}
\begin{split}
H^i(G_r,M)=&\varinjlim_{(n,p)=1}H^i(G_r,M[n]) \\
\cong&\varinjlim_{(n,p)=1}H^i(G_r,\Z/n\Z)\otimes_{\Z/n\Z}M[n] \\
\cong&\varinjlim_{(n,p)=1}(\bigwedge^iH^1(G_r,\Z/n\Z))\otimes_{\Z/n\Z}M[n]. \\
\end{split}
\end{equation}
We fix an isomorphism $G_r=\mr{Hom}(X_r,\hat{\Z}')$ with $X_r:=\Z^r$, and thus 
\[H^1(G_r,\Z/n\Z)=\mr{Hom}(G_r,\Z/n\Z)=X_r\otimes_{\Z}\Z/n\Z.\] 
Therefore we further have
\begin{align*}
H^i(G_r,M)\cong &\varinjlim_{(n,p)=1}(\bigwedge^i(X_r\otimes_{\Z}\Z/n\Z)) \otimes_{\Z/n\Z}M[n] \\
=&\varinjlim_{(n,p)=1}(\bigwedge^iX_r)\otimes_{\Z}M[n]  \\
=&(\bigwedge^iX_r)\otimes_{\Z}M',
\end{align*}
where the last two wedges are for $\Z$-module structure. Apparently the identification (\ref{eqC.2}) is induced by cup-product. 

To sum up, we have the following lemma.

\begin{lem}\label{C.1}
Let $p$ be a fixed prime number, and $\hat{\Z}':=\varprojlim_{(p,m)=1}\Z/m\Z$. Let $M$ be a torsion abelian group, and we regard it as a $(\hat{\Z}')^r$-module with respect to the trivial action. Let $M[p^\infty]$ (resp. $M'$, resp. $M[n]$) denote the $p$-primary part (resp. prime to $p$ part, resp. $n$-torsion part) of $M$. We fix an isomorphism $(\hat{\Z}')^r=\mr{Hom}(X_r,\hat{\Z}')$ with $X_r:=\Z^r$. 

Then 
\begin{enumerate}
\item for any positive integer $n$ with $(n,p)=1$, we have
\[H^1((\hat{\Z}')^r,M[n])=\mr{Hom}((\hat{\Z}')^r,M[n])\cong M[n]\otimes_{\Z}X_r\]
and the cup product induces an isomorphism
\[H^i((\hat{\Z}')^r,M[n])\cong M[n]\otimes_{\Z/n\Z}\bigwedge^iH^1((\hat{\Z}')^r,\Z/n\Z),\]
and the latter can be further identified with $M[n]\otimes_{\Z}(\bigwedge^iX_r) $;
\item for the profinite group cohomology of $(\hat{\Z}')^r$ with coefficients in $M$, we have
\begin{align*}
H^i((\hat{\Z}')^r,M)=\varinjlim_{(n,p)=1}H^i((\hat{\Z}')^r,M[n])=&\varinjlim_{(n,p)=1}M[n]\otimes_{\Z}(\bigwedge^iX_r) \\
=&M'\otimes_{\Z}(\bigwedge^iX_r).
\end{align*}
\end{enumerate}
\end{lem}

\begin{rmk}
Clearly the above computation works also for the profinite groups $\hat{\Z}^r$ and $\Z_l^r$ for any prime number $l$. Such results are the profinite group cohomology analogues of the description of the singular cohomology of topological tori (see \cite[Chap. 3, Exa. 3.16]{hat1}).
\end{rmk}

\section*{Acknowledgement}
The author thanks Professor Chicara Nakayama for very helpful discussions. This work was partially supported by the Research Training Group 2553 of the German Research Foundation DFG.

\bibliographystyle{alpha}
\bibliography{bib}

\end{document}